\documentclass[10pt]{amsart}
\usepackage{amsmath}
\usepackage{amsxtra}
\usepackage{amscd}
\usepackage{amsthm}
\usepackage{amsfonts}
\usepackage{amssymb}
\usepackage{eucal}

\newtheorem{cor}[subsection]{Corollary}

\newtheorem{lem}[subsection]{Lemma}
\newtheorem{prop}[subsection]{Proposition}

\newtheorem{conj}[subsection]{Conjecture}
\newtheorem{thm}[subsection]{Theorem}

\theoremstyle{definition}

\theoremstyle{remark}

\newcommand{\nc}{\newcommand}
\nc{\renc}{\renewcommand} \nc{\ssec}{\subsection}
\nc{\sssec}{\subsubsection} \nc{\on}{\operatorname}

\nc\ol{\overline} \nc\ul{\underline} \nc\wt{\widetilde}
\nc\tboxtimes{\wt{\boxtimes}} \nc{\alp}{\alpha}

\nc{\ZZ}{{\mathbb Z}} \nc{\NN}{{\mathbb N}} \nc{\CC}{{\mathbb C}}
\nc{\OO}{{\mathbb O}} \renc{\SS}{{\mathbb S}} \nc{\DD}{{\mathbb
D}}

\nc{\Fq}{{\mathbb F}_q} \nc{\Fqb}{\ol{{\mathbb F}_q}}
\nc{\Ql}{\ol{{\mathbb Q}_\ell}} \nc{\id}{\text{id}} \nc\X{\mathcal
X}

\nc{\Hom}{\on{Hom}} \nc{\Lie}{\on{Lie}} \nc{\Loc}{\on{Loc}}
\nc{\Pic}{\on{Pic}} \nc{\Bun}{\on{Bun}} \nc{\IC}{\on{IC}}
\nc{\Aut}{\on{Aut}} \nc{\rk}{\on{rk}} \nc{\Sh}{\on{Sh}}
\nc{\Perv}{\on{Perv}} \nc{\pos}{{\on{pos}}} \nc{\Conv}{\on{Conv}}
\nc{\Sph}{\on{Sph}} \nc{\Sym}{\on{Sym}}
\nc{\BunBb}{\overline{\Bun}_B} \nc{\Buno}{\overset{o}{\Bun}}
\nc{\BunPb}{{\overline{\Bun}_P}}
\nc{\BunBM}{\overline{\Bun}_{B(M)}}
\nc{\BunPbw}{{\widetilde{\Bun}_P}}
\nc{\BunBP}{\widetilde{\Bun}_{B,P}} \nc{\GUb}{\overline{G/U}}
\nc{\GUPb}{\overline{G/U(P)}}

\nc{\Hhom}{\underline{\on{Hom}}} \nc\syminfty{\on{Sym}^{\infty}}
\nc\lal{\ol{\lambda}} \nc\xl{\ol{x}} \nc\thl{\ol{\theta}}
\nc\nul{\ol{\nu}} \nc\mul{\ol{\mu}} \nc\Sum\Sigma
\nc{\oX}{\overset{o}{X}{}}

\nc{\M}{{\mathcal M}} \nc{\N}{{\mathcal N}} \nc{\F}{{\mathcal F}}
\nc{\D}{{\mathcal D}} \nc{\Q}{{\mathcal Q}} \nc{\Y}{{\mathcal Y}}
\nc{\G}{{\mathcal G}} \nc{\E}{{\mathcal E}} \nc{\CalC}{{\mathcal
C}}
\nc\Dh{\widehat{\D}}

\nc{\C}{{\mathcal C}} \nc{\K}{{\mathcal K}}
\renewcommand{\H}{{\mathcal H}}

\nc{\T}{{\mathcal T}} \nc{\V}{{\mathcal V}} \renc{\P}{{\mathcal
P}} \nc{\A}{{\mathcal A}} \nc{\B}{{\mathcal B}} \nc{\U}{{\mathcal
U}}

\nc{\Gr}{\on{Gr}}

\nc{\frn}{{\check{\mathfrak u}(P)}}
\nc\f{{\mathfrak f}}

\nc{\q}{{\mathfrak q}} \nc{\p}{{\mathfrak p}} \nc{\s}{{\mathfrak
s}} \nc\w{\text{w}}

\nc\mathi\iota \nc\Spec{\on{Spec}} \nc\Mod{\on{Mod}}
\nc{\tw}{\widetilde{\mathfrak t}} \nc{\pw}{\widetilde{\mathfrak
p}} \nc{\qw}{\widetilde{\mathfrak q}} \nc{\jw}{\widetilde j}

\nc{\grb}{\overline{\Gr}} \nc{\I}{\mathcal I}

\nc{\lambdach}{{\check\lambda}} \nc{\Lambdach}{{\check\Lambda}{}}
\nc{\much}{{\check\mu}} \nc{\omegach}{{\check\omega}}
\nc{\nuch}{{\check\nu}} \nc{\etach}{{\check\eta}}
\nc{\alphach}{{\check\alpha}} \nc{\betach}{{\check\beta}}
\nc{\rhoch}{{\check\rho}} \nc{\ch}{{\check h}}

\nc{\Hb}{\overline{\H}}

\nc{\BA}{{\mathbb{A}}} \nc{\BC}{{\mathbb{C}}}
\nc{\BM}{{\mathbb{M}}} \nc{\BN}{{\mathbb{N}}}
\nc{\BP}{{\mathbb{P}}} \nc{\BR}{{\mathbb{R}}}
\nc{\BZ}{{\mathbb{Z}}} \nc{\BS}{{\mathbb{S}}}

\nc{\CA}{{\mathcal{A}}} \nc{\CB}{{\mathcal{B}}}
\nc{\CE}{{\mathcal{E}}} \nc{\CF}{{\mathcal{F}}}
\nc{\CG}{{\mathcal{G}}} \nc{\CH}{{\mathcal{H}}}
\nc{\CI}{{\mathcal{I}}} \nc{\CL}{{\mathcal{L}}}
\nc{\CM}{{\mathcal{M}}} \nc{\CN}{{\mathcal{N}}}
\nc{\CO}{{\mathcal{O}}} \nc{\CP}{{\mathcal{P}}}
\nc{\CQ}{{\mathcal{Q}}} \nc{\CR}{{\mathcal{R}}}
\nc{\CS}{{\mathcal{S}}} \nc{\CT}{{\mathcal{T}}}
\nc{\CU}{{\mathcal{U}}} \nc{\CV}{{\mathcal{V}}}
\nc{\CW}{{\mathcal{W}}} \nc{\CZ}{{\mathcal{Z}}}

\nc{\cM}{{\check{\mathcal M}}{}} \nc{\csM}{{\check{\mathcal A}}{}}
\nc{\oM}{{\overset{\circ}{\mathcal M}}{}}
\nc{\obM}{{\overset{\circ}{\mathbf M}}{}}
\nc{\oCA}{{\overset{\circ}{\mathcal A}}{}}
\nc{\obA}{{\overset{\circ}{\mathbf A}}{}}
\nc{\ooM}{{\overset{\circ}{M}}{}}
\nc{\osM}{{\overset{\circ}{\mathsf M}}{}}
\nc{\vM}{{\overset{\bullet}{\mathcal M}}{}}
\nc{\nM}{{\underset{\bullet}{\mathcal M}}{}}
\nc{\oD}{{\overset{\circ}{\mathcal D}}{}}
\nc{\obD}{{\overset{\circ}{\mathbf D}}{}}
\nc{\oA}{{\overset{\circ}{\mathbb A}}{}}
\nc{\op}{{\overset{\bullet}{\mathbf p}}{}}
\nc{\cp}{{\overset{\circ}{\mathbf p}}{}}
\nc{\oU}{{\overset{\bullet}{\mathcal U}}{}}
\nc{\oZ}{{\overset{\circ}{\mathcal Z}}{}}
\nc{\ofZ}{{\overset{\circ}{\mathfrak Z}}{}}

\nc{\fa}{{\mathfrak{a}}} \nc{\fb}{{\mathfrak{b}}}
\nc{\fd}{{\mathfrak{d}}} \nc{\fe}{{\mathfrak{e}}}
\nc{\fg}{{\mathfrak{g}}} \nc{\fgl}{{\mathfrak{gl}}}
\nc{\fh}{{\mathfrak{h}}} \nc{\fri}{{\mathfrak{i}}}
\nc{\fj}{{\mathfrak{j}}} \nc{\fk}{{\mathfrak{k}}}
\nc{\fm}{{\mathfrak{m}}} \nc{\fn}{{\mathfrak{n}}}
\nc{\ft}{{\mathfrak{t}}} \nc{\fu}{{\mathfrak{u}}}
\nc{\fw}{{\mathfrak{w}}} \nc{\fz}{{\mathfrak{z}}}
\nc{\fp}{{\mathfrak{p}}} \nc{\frr}{{\mathfrak{r}}}
\nc{\fs}{{\mathfrak{s}}} \nc{\fsl}{{\mathfrak{sl}}}
\nc{\hsl}{{\widehat{\mathfrak{sl}}}}
\nc{\hgl}{{\widehat{\mathfrak{gl}}}}
\nc{\hg}{{\widehat{\mathfrak{g}}}}
\nc{\chg}{{\widehat{\mathfrak{g}}}{}^\vee}
\nc{\hn}{{\widehat{\mathfrak{n}}}}
\nc{\chn}{{\widehat{\mathfrak{n}}}{}^\vee}

\nc{\fA}{{\mathfrak{A}}} \nc{\fB}{{\mathfrak{B}}}
\nc{\fD}{{\mathfrak{D}}} \nc{\fE}{{\mathfrak{E}}}
\nc{\fF}{{\mathfrak{F}}} \nc{\fG}{{\mathfrak{G}}}
\nc{\fI}{{\mathfrak{I}}} \nc{\fJ}{{\mathfrak{J}}}
\nc{\fK}{{\mathfrak{K}}} \nc{\fL}{{\mathfrak{L}}}
\nc{\fM}{{\mathfrak{M}}} \nc{\fN}{{\mathfrak{N}}}
\nc{\frP}{{\mathfrak{P}}} \nc{\fQ}{{\mathfrak{Q}}}
\nc{\fT}{{\mathfrak{T}}} \nc{\fU}{{\mathfrak{U}}}
\nc{\fV}{{\mathfrak{V}}} \nc{\fW}{{\mathfrak{W}}}
\nc{\fX}{{\mathfrak{X}}} \nc{\fY}{{\mathfrak{Y}}}
\nc{\fZ}{{\mathfrak{Z}}}

\nc{\bb}{{\mathbf{b}}} \nc{\bc}{{\mathbf{c}}}
\nc{\be}{{\mathbf{e}}} \nc{\bj}{{\mathbf{j}}}
\nc{\bn}{{\mathbf{n}}} \nc{\bp}{{\mathbf{p}}}
\nc{\bq}{{\mathbf{q}}} \nc{\br}{{\mathbf{r}}}
\nc{\bfu}{{\mathbf{u}}} \nc{\bv}{{\mathbf{v}}}
\nc{\bx}{{\mathbf{x}}} \nc{\by}{{\mathbf{y}}}
\nc{\bw}{{\mathbf{w}}} \nc{\bA}{{\mathbf{A}}}
\nc{\bB}{{\mathbf{B}}} \nc{\bC}{{\mathbf{C}}}
\nc{\bD}{{\mathbf{D}}} \nc{\bF}{{\mathbf{F}}}
\nc{\bH}{{\mathbf{H}}} \nc{\bK}{{\mathbf{K}}}
\nc{\bM}{{\mathbf{M}}} \nc{\bN}{{\mathbf{N}}}
\nc{\bO}{{\mathbf{O}}} \nc{\bS}{{\mathbf{S}}}
\nc{\bV}{{\mathbf{V}}} \nc{\bW}{{\mathbf{W}}}
\nc{\bX}{{\mathbf{X}}} \nc{\bP}{{\mathbf{P}}}
\nc{\bZ}{{\mathbf{Z}}}

\nc{\sA}{{\mathsf{A}}} \nc{\sB}{{\mathsf{B}}}
\nc{\sC}{{\mathsf{C}}} \nc{\sD}{{\mathsf{D}}}
\nc{\sE}{{\mathsf{E}}} \nc{\sF}{{\mathsf{F}}}
\nc{\sK}{{\mathsf{K}}} \nc{\sL}{{\mathsf{L}}}
\nc{\sM}{{\mathsf{M}}} \nc{\sO}{{\mathsf{O}}}
\nc{\sQ}{{\mathsf{Q}}} \nc{\sP}{{\mathsf{P}}}
\nc{\sT}{{\mathsf{T}}} \nc{\sZ}{{\mathsf{Z}}}
\nc{\sfp}{{\mathsf{p}}} \nc{\sr}{{\mathsf{r}}}
\nc{\st}{{\mathsf{t}}} \nc{\sfb}{{\mathsf{b}}}
\nc{\sfc}{{\mathsf{c}}} \nc{\sd}{{\mathsf{d}}}
\nc{\sz}{{\mathsf{z}}}

\nc{\BK}{{\bar{K}}}

\nc{\tA}{{\widetilde{\mathbf{A}}}}
\nc{\tB}{{\widetilde{\mathcal{B}}}}
\nc{\tg}{{\widetilde{\mathfrak{g}}}} \nc{\tG}{{\widetilde{G}}}
\nc{\TM}{{\widetilde{\mathbb{M}}}{}}
\nc{\tO}{{\widetilde{\mathsf{O}}}{}}
\nc{\tU}{{\widetilde{\mathfrak{U}}}{}} \nc{\TZ}{{\tilde{Z}}}
\nc{\tx}{{\tilde{x}}} \nc{\tbv}{{\tilde{\bv}}}
\nc{\tfP}{{\widetilde{\mathfrak{P}}}{}} \nc{\tz}{{\tilde{\zeta}}}
\nc{\tmu}{{\tilde{\mu}}}

\nc{\urho}{\underline{\rho}} \nc{\uB}{\underline{B}}
\nc{\uC}{{\underline{\mathbb{C}}}} \nc{\ui}{\underline{i}}
\nc{\uj}{\underline{j}} \nc{\ofP}{{\overline{\mathfrak{P}}}}
\nc{\oB}{{\overline{\mathcal{B}}}}
\nc{\og}{{\overline{\mathfrak{g}}}} \nc{\oI}{{\overline{I}}}

\nc{\eps}{\varepsilon} \nc{\hrho}{{\hat{\rho}}}

\nc{\one}{{\mathbf{1}}} \nc{\two}{{\mathbf{t}}}

\nc{\Rep}{{\mathop{\operatorname{\rm Rep}}}}
\nc{\Tot}{{\mathop{\operatorname{\rm Tot}}}}
\nc{\Ker}{{\mathop{\operatorname{\rm Ker}}}}
\nc{\Hilb}{{\mathop{\operatorname{\rm Hilb}}}}
\nc{\End}{{\mathop{\operatorname{\rm End}}}}
\nc{\Ext}{{\mathop{\operatorname{\rm Ext}}}}
\nc{\CHom}{{\mathop{\operatorname{{\mathcal{H}}\it om}}}}
\nc{\GL}{{\mathop{\operatorname{\rm GL}}}}
\nc{\gr}{{\mathop{\operatorname{\rm gr}}}}
\nc{\Id}{{\mathop{\operatorname{\rm Id}}}}
\nc{\defi}{{\mathop{\operatorname{\rm def}}}}
\nc{\length}{{\mathop{\operatorname{\rm length}}}}
\nc{\supp}{{\mathop{\operatorname{\rm supp}}}}

\nc{\Cliff}{{\mathsf{Cliff}}}
\nc{\Fl}{{\mathsf{Fl}}} \nc{\Fib}{{\mathsf{Fib}}}
\nc{\Coh}{{\mathsf{Coh}}} \nc{\FCoh}{{\mathsf{FCoh}}}

\nc{\reg}{{\text{\rm reg}}}

\nc{\cplus}{{\mathbf{C}_+}} \nc{\cminus}{{\mathbf{C}_-}}
\nc{\cthree}{{\mathbf{C}_*}} \nc{\Qbar}{{\bar{Q}}}

\nc{\bh}{{\bar{h}}} \nc{\bOmega}{{\overline{\Omega}}}

\nc{\seq}[1]{\stackrel{#1}{\sim}}

\nc{\aff}{\operatorname{aff}}

\begin{document}

\author{Alexander Braverman and Michael Finkelberg}
\title
{Finite difference quantum Toda lattice via equivariant
$K$-theory}

\dedicatory{To Vladimir Drinfeld with admiration}



\address{{\it Address}:\newline
A.B.: Dept. of Math., Brown Univ., Providence, RI 02912, Einstein
Institute of Mathematics Edmond J. Safra Campus, Givat Ram The
Hebrew University of Jerusalem Jerusalem, 91904, Israel \newline
M.F.: Independent Moscow Univ., 11 Bolshoj Vlasjevskij per.,
Moscow 119002, Russia}

\email{\newline braval@math.brown.edu; fnklberg@mccme.ru}

\begin{abstract}
We construct the action of the quantum group
$U_v(\mathfrak{sl}_n)$ by the natural correspondences in the
equivariant localized $K$-theory of the Laumon based Quasiflags'
moduli spaces. The resulting module is the universal Verma module.
We construct geometrically the Shapovalov scalar product and the
Whittaker vectors. It follows that a certain generating function
of the characters of the global sections of the structure sheaves
of the Laumon moduli spaces satisfies a $v$-difference analogue of
the quantum Toda lattice system, reproving the main theorem of
Givental-Lee (cf. \cite{gl}). Similar constructions are performed
for the affine Lie agebra $\widehat{\mathfrak{sl}}_n$.
\end{abstract}
\maketitle

\section{Introduction}

\subsection{}
This work arose from an attempt to understand the results of the
paper ~\cite{gl} of A.~Givental and Y.-P.~Lee where the authors
perform some computations related to ``quantum $K$-theory" of flag
varieties (as well as some results from ~\cite{neok} related to 5d
$SU(n)$-gauge theory compactified on a circle) in the framework of
representation theory. Similar approach to quantum cohomology of
flag varieties (and to partition functions of 4d gauge theory) is
discussed in ~\cite{b} and ~\cite{be}.

In \cite{gl} the authors consider the moduli spaces $\fQ_{\ul{d}}$
introduced by G.~Laumon in ~\cite{la1}, ~\cite{la2}. These are
certain closures of the moduli spaces of based maps of degree
$\ul{d}$ from $\BP^1$ to the flag variety $\CB$ of
$\mathfrak{sl}_n$.

A Cartan torus $T$ of $SL_n$ acts on $\fQ_{\ul{d}}$. The
multiplicative group $\BC^*$ of dilations of $\BP^1$ (loop
rotations) also acts on $\fQ_{\ul{d}}$. The formal character of
the (infinite dimensional) $T\times\BC^*$-module
$R\Gamma(\fQ_{\ul{d}},\CO_{\ul{d}})$ turns out to be a rational
function on $T\times\BC^*$.
One may form a certain generating function $\fJ$ of these rational
functions for all degrees $\ul{d}$. Computing the function $\fJ$
presumably should give rise to a computation of the
$SL_n$-equivariant quantum $K$-theory ring of $\CB$ (which to the
best of the authors' knowledge has not yet been defined in the
literature).

A.~Givental and Y.-P.~Lee prove that $\fJ$ satisfies a certain
$v$-difference version of the quantum Toda lattice equations (here
$v$ stands for the tautological character of $\BC^*$). Moreover,
they suggest another way to construct solutions of the
$v$-difference Toda system: as the Shapovalov scalar product of
the Whittaker vectors in the universal Verma module for the
quantum group $U_v(\mathfrak{sl}_n)$. The latter construction was
worked out independently in ~\cite{e}, ~\cite{s2}.

\subsection{}
The principal goal of the present paper is to identify these two
constructions of solutions of the $v$-difference Toda system.
Namely, we prove that the natural correspondences between the
moduli spaces $\fQ_{\ul{d}}$ (for the degrees differing by a
simple root) give rise to the action of the standard generators of
$U_v(\mathfrak{sl}_n)$ on the localized equivariant $K$-theory
$\oplus_{\ul{d}}\ul{K}^{T\times\BC^*} (\fQ_{\ul{d}})$. Here the
localization is taken with respect to the
$K^{T\times\BC^*}(\cdot)=\BC[T\times\BC^*]$, that is, we tensor
everything with the fraction field of $\BC[T\times\BC^*]$. This is
needed since the above correspondences are not proper, but the
subspaces of their $T\times\BC^*$-fixed points are proper (in
fact, they are finite), so their action is well defined only in
the localized equivariant $K$-theory. This way we get a
$U_v(\mathfrak{sl}_n)$-module, and we identify it with the
universal Verma module $M$. We also compute in geometric terms the
Shapovalov scalar product on $M$, and the Whittaker vectors. It
turns out that the generating function for the Shapovalov scalar
product of the Whittaker vectors is a simple modification of the
Givental-Lee generating function $\fJ$. Thus we reprove the Main
Theorem of Givental-Lee.

\subsection{}
There is a similar generating function $J$ for equivariant
integrals of the unit cohomology classes of $\fQ_{\ul{d}}$ which
controls the $T$-equivariant quantum cohomology of $\CB$. It
satisfies the quantum Toda lattice differential system, as proved
originally by A.~Givental and B.~Kim.

For the simple Lie algebras $\fg$ other than $\mathfrak{sl}_n$
there is no analogue of the Laumon moduli spaces $\fQ_{\ul{d}}$
but there is Drinfeld's moduli space of Quasimaps
$\CZ_{\ul{d}}(\fg)$. It also exists for the case of affine Lie
algebras, under the name of Uhlenbeck compactification. In the
$\mathfrak{sl}_n$ case $\fQ_{\ul{d}}$ is a small resolution of
$\CZ_{\ul{d}}(\mathfrak{sl}_n)$. In the affine
$\widehat{\mathfrak{sl}}_n$ case
$\CZ_{\ul{d}}(\widehat{\mathfrak{sl}}_n)$ possesses a semismall
resolution of singularities: the moduli space $\CP_{\ul{d}}$ of
torsion free parabolic sheaves on $\BP^1\times\BP^1$ endowed with
some additional structures. Thus in the affine case we can define
an analog of the function $J$ which we denote by $J_{\aff}$ (this
is discussed in \cite{b}).

The generating function $J$ (for any simple $G$) is known to
satisfy the quantum (differential) Toda equations (cf. \cite{gk}
and \cite{kim}).

In the work ~\cite{b}, the generating function $J$ for equivariant
integrals of the unit cohomology classes of $\CZ_{\ul{d}}(\fg)$
was proved to satisfy the quantum Toda lattice by constructing the
action of the Langlands dual Lie algebra $\check\fg$ in the
equivariant Intersection Cohomology of the Drinfeld
compactifications. Also in the affine case the function $J_{\aff}$
was shown to satisfy some non-stationary analog of ``the most
basic" (quadratic) Toda equation. Thus ~\cite{b} offered a
representation theoretic explanation of the Givental-Kim results
as well as generalized them to the affine case. And the present
work is a multiplicative analogue of ~\cite{b} in the simplest
case of $\mathfrak{sl}_n$.

\subsection{}
It would be extremely interesting to extend our work to other
simple and affine Lie algebras. It would require something like an
equivariant ``IC $K$-theory'' of $\CZ_{\ul{d}}(\fg)$ which is not
defined at the moment. In case of $\mathfrak{sl}_n$ the IC
cohomology of $\CZ_{\ul{d}}(\mathfrak{sl}_n)$ coincides with the
cohomology of the small resolution $\fQ_{\ul{d}}$, while in the
affine case of $\widehat{\mathfrak{sl}}_n$ the IC cohomology of
$\CZ_{\ul{d}}(\widehat{\mathfrak{sl}}_n)$ is a direct summand in
the cohomology of the semismall resolution $\CP_{\ul{d}}$.
Accordingly, one might look for the correct ``IC $K$-theory'' of
$\CZ_{\ul{d}}(\widehat{\mathfrak{sl}}_n)$ as an appropriate direct
summand of the usual $K$-theory of $\CP_{\ul{d}}$.

This is sketched in the Section ~\ref{p}. Namely, similarly to the
case of Laumon spaces, the quantum affine group
$U_v(\widehat{\mathfrak{sl}}_n)$ acts by the natural
correspondences on the direct sum of localized equivariant
$K$-groups $\oplus_{\ul{d}}\ul{K}^{T\times\BC^*\times\BC^*}
(\CP_{\ul{d}})$. However, this module looks more like the
universal Verma module for $U_v(\widehat{\mathfrak{gl}}_n)$, and
we have to specify a certain submodule isomorphic to the universal
Verma module for $U_v(\widehat{\mathfrak{sl}}_n)$. Then we
construct geometrically the Shapovalov scalar product, and the
Whittaker vectors. It turns out that the Shapovalov scalar product
of the Whittaker vectors can be expressed via the formal
characters of the global sections
$R\Gamma(\CP_{\ul{d}},\CO_{\ul{d}})$ as in the case of
$\mathfrak{sl}_n$. However, we were unable to derive any
$v$-difference equation for the affine version of the generating
function $\fJ$.

\subsection{Acknowledgments}
M.F. is obliged to V.~Schechtman, A.~Stoyanovsky, B.~Feigin,
E.~Vasserot, and R.~Bezrukavnikov who, ever since the appearance
of ~\cite{fk}, urged him to consider its equivariant $K$-theory
analogue. While trying to guess the correct formulae in the low
ranks, we profited strongly from the computational help of
V.~Dotsenko, V.~Golyshev, A.~Kuznetsov. We are also grateful to
P.~Etingof and A.~Joseph for very useful explanations; to
M.~Kashiwara for bringing the reference ~\cite{nz} to our
attention, and to the referee for the valuable comments.
Last but not least, our thanks go to A.~Tsymbaliuk for the careful
reading of our note and spotting several mistakes.
We would like to thank the Weizmann Institute and RIMS, Kyoto, as well as
the University of Chicago, for the hospitality and support.

M.F. was partially supported by the CRDF award RM1-2545-MO-03.
A.B. was partially supported by the NSF grant DMS-0300271.

\section{Laumon spaces and quantum groups}

\subsection{}
We recall the setup of ~\cite{fk}. Let $\bC$ be a smooth
projective curve of genus zero. We fix a coordinate $z$ on $\bC$,
and consider the action of $\BC^*$ on $\bC$ such that
$v(z)=v^{-2}z$. We have $\bC^{\BC^*}=\{0,\infty\}$.

We consider an $n$-dimensional vector space $W$ with a basis
$w_1,\ldots,w_n$. This defines a Cartan torus $T\subset
G=SL_n\subset Aut(W)$. We also consider its $2^{n-1}$-fold cover,
the bigger torus $\widetilde{T}$, acting on $W$ as follows: for
$\widetilde{T}\ni\ul{t}=(t_1,\ldots,t_n)$ we have
$\ul{t}(w_i)=t_i^2w_i$. We denote by $\CB$ the flag variety of
$G$.

\subsection{}
Given an $(n-1)$-tuple of nonnegative integers
$\ul{d}=(d_1,\ldots,d_{n-1})$, we consider the Laumon's
quasiflags' space $\CQ_{\ul{d}}$, see ~\cite{la2}, ~4.2. It is the
moduli space of flags of locally free subsheaves
$$0\subset\CW_1\subset\ldots\subset\CW_{n-1}\subset\CW=W\otimes\CO_\bC$$
such that $\on{rank}(\CW_k)=k$, and $\deg(\CW_k)=-d_k$.

It is known to be a smooth projective variety of dimension
$2d_1+\ldots+2d_{n-1}+\dim\CB$, see ~\cite{la1}, ~2.10.

\subsection{} We consider the following locally closed subvariety
$\fQ_{\ul{d}}\subset\CQ_{\ul{d}}$ (quasiflags based at
$\infty\in\bC$) formed by the flags
$$0\subset\CW_1\subset\ldots\subset\CW_{n-1}\subset\CW=W\otimes\CO_\bC$$
such that $\CW_i\subset\CW$ is a vector subbundle in a
neighbourhood of $\infty\in\bC$, and the fiber of $\CW_i$ at
$\infty$ equals the span $\langle w_1,\ldots,w_i\rangle\subset W$.

It is known to be a smooth quasiprojective variety of dimension
$2d_1+\ldots+2d_{n-1}$.

\subsection{}
\label{fixed points} The group $G\times\BC^*$ acts naturally on
$\CQ_{\ul{d}}$, and the group $\widetilde{T}\times\BC^*$ acts
naturally on $\fQ_{\ul{d}}$. The set of fixed points of
$\widetilde{T}\times\BC^*$ on $\fQ_{\ul{d}}$ is finite; we recall
its description from ~\cite{fk}, ~2.11.

Let $\widetilde{\ul{d}}$ be a collection of nonnegative integers
$(d_{ij}),\ i\geq j$, such that $d_i=\sum_{j=1}^id_{ij}$, and for
$i\geq k\geq j$ we have $d_{kj}\geq d_{ij}$. Abusing notation we
denote by $\widetilde{\ul{d}}$ the corresponding
$\widetilde{T}\times\BC^*$-fixed point in $\fQ_{\ul{d}}$:

$\CW_1=\CO_\bC(-d_{11}\cdot0)w_1,$

$\CW_2=\CO_\bC(-d_{21}\cdot0)w_1\oplus\CO_\bC(-d_{22}\cdot0)w_2,$

$\ldots\ \ldots\ \ldots\ ,$

$\CW_{n-1}=\CO_\bC(-d_{n-1,1}\cdot0)w_1\oplus\CO_\bC(-d_{n-1,2}\cdot0)w_2
\oplus\ldots\oplus\CO_\bC(-d_{n-1,n-1}\cdot0)w_{n-1}.$

\subsection{}
For $i\in\{1,\ldots,n-1\}$, and $\ul{d}=(d_1,\ldots,d_{n-1})$, we
set $\ul{d}+i:=(d_1,\ldots,d_i+1,\ldots,d_{n-1})$. We have a
correspondence $\CE_{\ul{d},i}\subset\CQ_{\ul{d}}\times
\CQ_{\ul{d}+i}$ formed by the pairs $(\CW_\bullet,\CW'_\bullet)$
such that for $j\ne i$ we have $\CW_j=\CW'_j$, and
$\CW'_i\subset\CW_i$, see ~\cite{fk}, ~3.1. In other words,
$\CE_{\ul{d},i}$ is the moduli space of flags of locally free
sheaves
$$0\subset\CW_1\subset\ldots\CW_{i-1}\subset\CW'_i\subset\CW_i\subset
\CW_{i+1}\ldots\subset\CW_{n-1}\subset\CW$$ such that
$\on{rank}(\CW_k)=k$, and $\deg(\CW_k)=-d_k$, while
$\on{rank}(\CW'_i)=i$, and $\deg(\CW'_i)=-d_i-1$.

According to ~\cite{la1}, ~2.10, $\CE_{\ul{d},i}$ is a smooth
projective algebraic variety of dimension
$2d_1+\ldots+2d_{n-1}+\dim\CB+1$.

We denote by $\bp$ (resp. $\bq$) the natural projection
$\CE_{\ul{d},i}\to\CQ_{\ul{d}}$ (resp.
$\CE_{\ul{d},i}\to\CQ_{\ul{d}+i}$). We also have a map $\br:\
\CE_{\ul{d},i}\to\bC,$
$$(0\subset\CW_1\subset\ldots\subset\CW_{i-1}\subset\CW'_i\subset\CW_i\subset
\CW_{i+1}\subset\ldots\subset\CW_{n-1}\subset\CW)\mapsto\on{supp}(\CW_i/\CW'_i).$$

The correspondence $\CE_{\ul{d},i}$ comes equipped with a natural
line bundle $\CL_i$ whose fiber at a point
$$(0\subset\CW_1\subset\ldots\subset\CW_{i-1}\subset\CW'_i\subset\CW_i\subset
\CW_{i+1}\subset\ldots\subset\CW_{n-1}\subset\CW)$$ equals
$\Gamma(\bC,\CW_i/\CW'_i)$.

Finally, we have a transposed correspondence
$^\sT\CE_{\ul{d},i}\subset \CQ_{\ul{d}+i}\times\CQ_{\ul{d}}$.

\subsection{}
Restricting to $\fQ_{\ul{d}}\subset\CQ_{\ul{d}}$ we obtain the
correspondence
$\fE_{\ul{d},i}\subset\fQ_{\ul{d}}\times\fQ_{\ul{d}+i}$ together
with line bundle $\fL_i$ and the natural maps $\bp:\
\fE_{\ul{d},i}\to\fQ_{\ul{d}},\ \bq:\
\fE_{\ul{d},i}\to\fQ_{\ul{d}+i},\ \br:\
\fE_{\ul{d},i}\to\bC-\infty$. We also have a transposed
correspondence $^\sT\fE_{\ul{d},i}\subset
\fQ_{\ul{d}+i}\times\fQ_{\ul{d}}$. It is a smooth quasiprojective
variety of dimension $2d_1+\ldots+2d_{n-1}+1$.

\subsection{}
We denote by ${}'M$ the direct sum of equivariant (complexified)
$K$-groups:
${}'M=\oplus_{\ul{d}}K^{\widetilde{T}\times\BC^*}(\fQ_{\ul{d}})$.
It is a module over
$K^{\widetilde{T}\times\BC^*}(pt)=\BC[\widetilde{T}\times\BC^*]=
\BC[t_1,\ldots,t_n,v\ :\ t_1\cdots t_n=1]$.
We define $M=\ {}'M\otimes_{K^{\widetilde{T}\times\BC^*}(pt)}
\on{Frac}(K^{\widetilde{T}\times\BC^*}(pt))$.

We have an evident grading $M=\oplus_{\ul{d}}M_{\ul{d}},\
M_{\ul{d}}=K^{\widetilde{T}\times\BC^*}(\fQ_{\ul{d}})
\otimes_{K^{\widetilde{T}\times\BC^*}(pt)}
\on{Frac}(K^{\widetilde{T}\times\BC^*}(pt))$.

\subsection{}
\label{operators} The grading and the correspondences
$^\sT\fE_{\ul{d},i},\fE_{\ul{d},i}$ give rise to the following
operators on $M$ (note that though $\bp$ is not proper, $\bp_*$ is
well defined on the localized equivariant $K$-theory due to the
finiteness of the fixed point sets):



$K_i=t_{i+1}t_i^{-1}v^{2d_i-d_{i-1}-d_{i+1}+1}:\ M_{\ul{d}}\to
M_{\ul{d}}$;



$L_i=t_1^{-1}\cdots t_i^{-1}v^{d_i+\frac{1}{2}i(n-i)}:\
M_{\ul{d}}\to M_{\ul{d}}$;

$f_i=\bp_*\bq^*:\ M_{\ul{d}}\to M_{\ul{d}-i}$;

$F_i=t_{i+1}^it_i^{-i}v^{2id_i-id_{i-1}-id_{i+1}-i}\bp_*\bq^*:\
M_{\ul{d}}\to M_{\ul{d}-i}$;

$e_i=-t_{i+1}^{-1}t_i^{-1}v^{d_{i+1}-d_{i-1}}\bq_*(\fL_i\otimes\bp^*):\
M_{\ul{d}}\to M_{\ul{d}+i}$,

$E_i=-t_{i+1}^{-i-1}t_i^{i-1}v^{(i-1)d_{i-1}+(i+1)d_{i+1}-2id_i-i}
\bq_*(\fL_i\otimes\bp^*):\ M_{\ul{d}}\to M_{\ul{d}+i}$.

\subsection{}
We recall the notations and results of ~\cite{s} in the special
case of quantum group of $SL_n$ type.

$U$ is the $\BC[v,v^{-1}]$-algebra with generators
$E_i,L_i^{\pm1},K_i^{\pm1},F_i,\ 1\leq i\leq n-1$, subject to the
following relations:

\begin{equation}
\label{och} L_iL_j=L_jL_i,\ K_1=L_1^2L_2^{-1},\
K_i=L_{i-1}^{-1}L_i^2L_{i+1}^{-1},\ K_{n-1}=L_{n-2}^{-1}L_{n-1}^2
\end{equation}

\begin{equation}
\label{ochev} L_iE_jL_i^{-1}=v^{\delta_{i,j}}E_j,\
L_iF_jL_i^{-1}=v^{-\delta_{i,j}}F_j
\end{equation}

\begin{equation}
\label{ochevidno}
E_iF_j-F_jE_i=\delta_{i,j}\frac{K_i-K_i^{-1}}{v-v^{-1}}
\end{equation}

\begin{equation}
\label{Serre1} |i-j|>1\ \Longrightarrow\
E_iE_j-E_jE_i=0=F_iF_j-F_jF_i
\end{equation}

\begin{equation}
\label{Serre2} |i-j|=1\ \Longrightarrow\
E_i^2E_j-(v+v^{-1})E_iE_jE_i+E_jE_i^2=0=
F_i^2F_j-(v+v^{-1})F_iF_jF_i+F_jF_i^2
\end{equation}

Sevostyanov considers elements $e_i,f_i\in U$ depending on a
choice of $(n-1)\times(n-1)$-matrices $n_{ij},\ c_{ij}$. We make
the following choice:

\begin{equation}
\label{vybor} n_{i,i}=-2i;\ n_{i,i+1}=n_{i,i-1}=i,
\end{equation}
otherwise $n_{ij}=0$.
\begin{equation}
\label{Rossii} i<n-1\ \Longrightarrow\ c_{i,i+1}=-1,\ c_{i+1,i}=1,
\end{equation}
otherwise $c_{ij}=0$. In other words, $c_{ij}=n_{ij}-n_{ji}$.

Then we have
\begin{equation}
\label{Seva} f_i:=L_{i-1}^iL_i^{-2i}L_{i+1}^iF_i=K_i^{-i}F_i,\
e_i:=E_iL_{i-1}^{-i}L_i^{2i}L_{i+1}^{-i}=E_iK_i^i.
\end{equation}

Clearly, the algebra $U$ is generated by
$e_i,L_i^{\pm1},K_i^{\pm1},f_i,\ 1\leq i\leq n-1$, and the
relations ~(\ref{ochev})--~(\ref{Serre2}) above are equivalent to
the relations ~(\ref{ochev'})--~(\ref{Serre2'}) below.

\begin{equation}
\label{ochev'} L_ie_jL_i^{-1}=v^{\delta_{i,j}}e_j,\
L_if_jL_i^{-1}=v^{-\delta_{i,j}}f_j
\end{equation}

\begin{equation}
\label{ochevidno'}
e_if_j-v^{c_{ij}}f_je_i=\delta_{i,j}\frac{K_i-K_i^{-1}}{v-v^{-1}}
\end{equation}

\begin{equation}
\label{Serre1'} |i-j|>1\ \Longrightarrow\
e_ie_j-e_je_i=0=f_if_j-f_jf_i
\end{equation}

\begin{equation}
\label{Serre2'} |i-j|=1\ \Longrightarrow\
e_i^2e_j-v^{c_{ij}}(v+v^{-1})e_ie_je_i+v^{2c_{ij}}e_je_i^2=0=
f_i^2f_j-v^{c_{ij}}(v+v^{-1})f_if_jf_i+v^{2c_{ij}}f_jf_i^2
\end{equation}

\subsection{Remark}
The elements $f_i$ of the subalgebra $U_{\leq0}$ generated by
$F_1,\ldots,F_{n-1},K_1,\ldots,K_{n-1}$ were introduced by
C.~M.~Ringel in ~\cite{r}. They are the natural generators of the
Hall algebra of the $A_{n-1}$-quiver with the set of vertices
$1,\ldots,n-1$, and orientation $i\longrightarrow i+1$. More
generally, Ringel's construction works for an arbitrary
orientation of an $ADE$ quiver, and produces Sevostyanov's
generators $f_i$ (in the simply laced case). It can be seen easily
that the set of Sevostyanov's matrices $c_{ij}$ (parametrizing the
choices of his ``Coxeter realizations'') is in a natural bijection
with the set of orientations of the corresponding quiver.

\subsection{}
We are finally able to formulate our main theorem. Recall the
operators $E_i,e_i,L_i^{\pm1},K_i^{\pm1},F_i,f_i$ on $M$ defined
in ~\ref{operators}.

\begin{thm}
\label{main} The operators $E_i,L_i^{\pm1},K_i^{\pm1},F_i,\ 1\leq
i\leq n-1$, on $M$ satisfy the relations
~(\ref{och})--(\ref{Serre2}). Equivalently, the operators
$e_i,L_i^{\pm1},K_i^{\pm1},f_i,\ 1\leq i\leq n-1$, on $M$ satisfy
the relations ~(\ref{och}), ~(\ref{ochev'})--(\ref{Serre2'}).
\end{thm}

The relations ~(\ref{och}) and ~(\ref{ochev}) are evident. The
relation ~(\ref{ochevidno}) for $i\ne j$ follows from a
transversality property formulated in the next subsection.

\subsection{}
We consider the subvarieties $\bp_{12}^{-1}(\fE_{\ul{d},i})$ and
$\bp_{23}^{-1}(\ ^\sT\fE_{\ul{d}+i-j,j})$ in
$\fQ_{\ul{d}}\times\fQ_{\ul{d}+i}\times\fQ_{\ul{d}+i-j}$.
Similarly, we consider the subvarieties $\bp_{12}^{-1}(\
^\sT\fE_{\ul{d}-j,j})$ and $\bp_{23}^{-1}(\fE_{\ul{d}-j,i})$ in
$\fQ_{\ul{d}}\times\fQ_{\ul{d}-j}\times\fQ_{\ul{d}+i-j}$.

\begin{lem}
\label{trans} For $i\ne j$ the intersection (a)
$\bp_{12}^{-1}(\fE_{\ul{d},i})\cap\bp_{23}^{-1}(\
^\sT\fE_{\ul{d}+i-j,j})$ in
$\fQ_{\ul{d}}\times\fQ_{\ul{d}+i}\times\fQ_{\ul{d}+i-j}$ (resp.
(b) $\bp_{12}^{-1}(\
^\sT\fE_{\ul{d}-j,j})\cap\bp_{23}^{-1}(\fE_{\ul{d}-j,i})$ in
$\fQ_{\ul{d}}\times\fQ_{\ul{d}-j}\times\fQ_{\ul{d}+i-j}$) is
transversal.

(c) $\bp_{12}^{-1}(\fE_{\ul{d},i})\cap\bp_{23}^{-1}(\
^\sT\fE_{\ul{d}+i-j,j}) \simeq \bp_{12}^{-1}(\
^\sT\fE_{\ul{d}-j,j})\cap\bp_{23}^{-1}(\fE_{\ul{d}-j,i})$.
\end{lem}

\begin{proof}
We prove (a). By definition, $\bp_{12}^{-1}(\fE_{\ul{d},i})$ is
the moduli space of pairs of flags
$$(0\subset\CW'_1=\CW_1\subset\CW'_2=\CW_2\subset
\ldots\subset\CW'_i\subset\CW_i
\subset\ldots\subset\CW'_{n-1}=\CW_{n-1}\subset\CW,$$
$$0\subset\CW'''_1\subset\CW'''_2\subset\ldots\subset\CW'''_{n-1}\subset\CW)$$
of prescribed ranks and degrees, while $\bp_{23}^{-1}(\
^\sT\fE_{\ul{d}+i-j,j})$ is the moduli space of pairs of flags
$$(0\subset\CW_1\subset\ldots\subset\CW_{n-1}\subset\CW,$$
$$0\subset\CW'_1=\CW'''_1\subset\CW'_2=\CW'''_2\subset
\ldots\subset\CW'_j\subset
\CW'''_j\subset\ldots\subset\CW'_{n-1}=\CW'''_{n-1}\subset\CW)$$
of prescribed ranks and degrees.

Their intersection is the moduli space of flags (say, $i<j$)
$$0\subset\CW'_1=\CW_1=\CW'''_1\subset\ldots\subset\CW'_i=\CW'''_i\subset\CW_i
\subset\ldots\subset\CW'_j=$$
$$=\CW_j\subset\CW'''_j\subset\ldots\subset
\CW'_{n-1}=\CW_{n-1}=\CW'''_{n-1}\subset\CW$$ of prescribed ranks
and degrees which is smooth according to ~\cite{la1}, ~2.10. This
implies that at any closed point of the scheme-theoretic
intersection $\bp_{12}^{-1}(\fE_{\ul{d},i})\cap\bp_{23}^{-1}(\
^\sT\fE_{\ul{d}+i-j,j})$ the Zariski tangent space to
$\bp_{12}^{-1}(\fE_{\ul{d},i})\cap\bp_{23}^{-1}(\
^\sT\fE_{\ul{d}+i-j,j})$ is the intersection of tangent spaces to
$\bp_{12}^{-1}(\fE_{\ul{d},i})$ and $\bp_{23}^{-1}(\
^\sT\fE_{\ul{d}+i-j,j})$. Comparing the dimensions we conclude
that the sum of tangent spaces to $\bp_{12}^{-1}(\fE_{\ul{d},i})$
and $\bp_{23}^{-1}(\ ^\sT\fE_{\ul{d}+i-j,j})$ must coincide with
the tangent space to
$\fQ_{\ul{d}}\times\fQ_{\ul{d}+i}\times\fQ_{\ul{d}+i-j}$. Hence
the intersection is transversal. This completes the proof of (a).

In (b) we prove similarly that $\bp_{12}^{-1}(\
^\sT\fE_{\ul{d}-j,j})\cap\bp_{23}^{-1}(\fE_{\ul{d}-j,i})$ is the
moduli space of flags (say, $i<j$)
$$0\subset\CW_1=\CW'''_1=\CW''_1\subset\ldots\subset\CW'''_i\subset\CW_i
=\CW''_i \subset\ldots\subset\CW_j\subset$$
$$\subset\CW'''_j=\CW''_j\subset \ldots\subset
\CW_{n-1}=\CW'''_{n-1}=\CW''_{n-1}\subset\CW$$ of prescribed ranks
and degrees which is smooth according to ~\cite{la1}, ~2.10. Hence
the intersection is transversal by the same argument as in the
proof of (a). This completes the proof of (b).

Part (c) was proved in ~\cite{fk}, ~3.6. We just recall that the
mutually inverse isomorphisms send a triple
$(\CW_\bullet,\CW'_\bullet,\CW'''_\bullet)$ to
$(\CW_\bullet,\CW''_\bullet,\CW'''_\bullet)$ where
$\CW''_\bullet:=\CW_\bullet+\CW'''_\bullet$, and a triple
$(\CW_\bullet,\CW''_\bullet,\CW'''_\bullet)$ to
$(\CW_\bullet,\CW'_\bullet,\CW'''_\bullet)$ where
$\CW'_\bullet:=\CW_\bullet\cap\CW'''_\bullet$.
\end{proof}

\subsection{}
We return to the proof of relation ~(\ref{ochevidno}) for $i\ne
j$. The composition $F_jE_i$ is given by the action of
correspondence
$$f(\ul{t})g(v)\bp_{13*}(\bp_{12}^*\fL_i
\stackrel{L}{\otimes}_{\CO_{\fQ_{\ul{d}}\times\fQ_{\ul{d}+i}\times
\fQ_{\ul{d}+i-j}}}\bp_{23}^*\CO_{^\sT\fE_{\ul{d}+i-j,j}})$$ where
$f$ (resp. $g$) is a certain monomial in $\ul{t}$ (resp. $v$).

Because of the transversality in ~\ref{trans}(a), $\bp_{12}^*\fL_i
\stackrel{L}{\otimes}_{\CO_{\fQ_{\ul{d}}\times\fQ_{\ul{d}+i}\times
\fQ_{\ul{d}+i-j}}}\bp_{23}^*\CO_{^\sT\fE_{\ul{d}+i-j,j}}$ is a
line bundle $\fL_{i,j}$ on
$\bp_{12}^{-1}(\fE_{\ul{d},i})\cap\bp_{23}^{-1}(\
^\sT\fE_{\ul{d}+i-j,j})$ whose fiber at a point
$(\CW_\bullet,\CW'_\bullet,\CW'''_\bullet)$ is equal to
$\Gamma(\bC,\CW_i/\CW'''_i)$.

Similarly, due to the transversality in ~\ref{trans}(b), the
composition $E_iF_j$ is given by the action of correspondence
$$f'(\ul{t})g'(v)\bp_{13*}(\fL'_{i,j})$$
where $f'$ (resp. $g'$) is a certain monomial in $\ul{t}$ (resp.
$v$), and $\fL'_{i,j}$ is a line bundle on $\bp_{12}^{-1}(\
^\sT\fE_{\ul{d}-j,j})\cap\bp_{23}^{-1}(\fE_{\ul{d}-j,i})$ whose
fiber at a point $(\CW_\bullet,\CW''_\bullet,\CW'''_\bullet)$ is
equal to $\Gamma(\bC,\CW_i/\CW'''_i)$.

Now the isomorphism in ~\ref{trans}(c) clearly takes $\fL_{i,j}$
to $\fL'_{i,j}$, and a routine check shows that
$f(\ul{t})g(v)=f'(\ul{t})g'(v)$. This completes the proof of the
relations ~(\ref{ochevidno}) for $i\ne j$.

\subsection{}
\label{matrix} To prove the relation ~(\ref{ochevidno}) for $i=j$
we use the localization to the fixed points.

According to the Thomason localization theorem (see e.g.
~\cite{cg}), restriction to the $\widetilde{T}\times\BC^*$-fixed
point set induces an isomorphism
$$K^{\widetilde{T}\times\BC^*}(\fQ_{\ul{d}})
\otimes_{K^{\widetilde{T}\times\BC^*}(pt)}
\on{Frac}(K^{\widetilde{T}\times\BC^*}(pt))\to
K^{\widetilde{T}\times\BC^*}(\fQ_{\ul{d}}^{\widetilde{T}\times\BC^*})
\otimes_{K^{\widetilde{T}\times\BC^*}(pt)}
\on{Frac}(K^{\widetilde{T}\times\BC^*}(pt))$$ (resp.
$$K^{\widetilde{T}\times\BC^*}(\fE_{\ul{d},i})
\otimes_{K^{\widetilde{T}\times\BC^*}(pt)}
\on{Frac}(K^{\widetilde{T}\times\BC^*}(pt))\to
K^{\widetilde{T}\times\BC^*}(\fE_{\ul{d},i}^{\widetilde{T}\times\BC^*})
\otimes_{K^{\widetilde{T}\times\BC^*}(pt)}
\on{Frac}(K^{\widetilde{T}\times\BC^*}(pt)))$$



The classes of the structure sheaves $[\widetilde{\ul{d}}]$ of the
$\widetilde{T}\times\BC^*$-fixed points $\widetilde{\ul{d}}$ (see
~\ref{fixed points}) form a basis in
$\oplus_{\ul{d}}K^{\widetilde{T}\times\BC^*}
(\fQ_{\ul{d}}^{\widetilde{T}\times\BC^*})
\otimes_{K^{\widetilde{T}\times\BC^*}(pt)}\on{Frac}
(K^{\widetilde{T}\times\BC^*}(pt))$. In order to compute the
matrix coefficients of $E_i,F_i$ in this basis, we have to know
the character of the $\widetilde{T}\times\BC^*$-action in the
tangent spaces $\CT_{\widetilde{\ul{d}}}\fQ_{\ul{d}}$ and also in
the tangent spaces to the fixed points in the correspondences.
This is the subject of the following Proposition.

\subsection{}

Note that a point $(\widetilde{\ul{d}},\widetilde{\ul{d}}{}')$
lies in the correspondence $\fE_{\ul{d},i}$ if and only if
$d_{k,j}=d'_{k,j}$ with a single exception $d'_{i,j}=d_{i,j}+1$
for certain $j\leq i$.

\begin{prop}
\label{zanudstvo} a) The character $\chi_{\widetilde{\ul{d}}}$ of
$\widetilde{T}\times\BC^*$ in the tangent space
$\CT_{\widetilde{\ul{d}}}\fQ_{\ul{d}}$ equals
$$\sum_{1\leq j<k\leq n}t_k^2t_j^{-2}\sum_{l=1}^{d_{k-1,j}-d_{k,k}}v^{2l}+
\sum_{1\leq j,k\leq n-1}t_k^2t_j^{-2}\sum_{i=\max(k,j)}^{n-1}
\sum_{l=d_{i,j}-d_{i,k}+1}^{d_{i,j}-d_{i+1,k}}v^{2l}$$
where we set $d_{n,k}=0$.

b) The character
$\chi_{(\widetilde{\ul{d}},\widetilde{\ul{d}}{}')}$ of
$\widetilde{T}\times\BC^*$ in the tangent space
$\CT_{(\widetilde{\ul{d}},\widetilde{\ul{d}}{}')}\fE_{\ul{d},i}$
equals
$$\chi_{\widetilde{\ul{d}}}+
\sum_{k\leq i}t_j^2t_k^{-2}v^{2d'_{i,k}-2d_{i,j}}- \sum_{k\leq
i-1}t_j^2t_k^{-2}v^{2d_{i-1,k}-2d_{i,j}}$$
if $d'_{i,j}=d_{i,j}+1$ for certain $j\leq i$.

c) The character
$\lambda_{(\widetilde{\ul{d}},\widetilde{\ul{d}}{}')}$ of
$\widetilde{T}\times\BC^*$ in the fiber of $\fL_i$ at the point
$(\widetilde{\ul{d}},\widetilde{\ul{d}}{}')$ equals
$t_j^2v^{-2d_{i,j}}$ if $d'_{i,j}=d_{i,j}+1$.
\end{prop}

\begin{proof}
Let $\CQ$ be the moduli space of flags of locally free subsheaves
$$0\subset\CW_1\subset\CW_2\subset\ldots\subset\CW_r\subset\CW$$ of fixed
ranks. Then the tangent space $\CT_{\CW_\bullet}\CQ$ equals the
kernel of
$$\sum_{1\leq l<r}p_{l-1}^*\otimes\on{Id}-\on{Id}\otimes q_l:\
\oplus_l\on{Hom}(\CW_l,\CW/\CW_l)\twoheadrightarrow
\oplus_l\on{Hom}(\CW_l,\CW/\CW_{l+1})$$ where $p_l:\
\CW_l\hookrightarrow\CW_{l+1};\ q_l:\
\CW/\CW_l\twoheadrightarrow\CW/\CW_{l+1}$ (see e.g. ~\cite{gl},
~3.2).

Now the parts a), b) follow easily from the obvious equalities
$\on{ch}(\on{Hom}(\CO_\bC(-a),\CO_\bC)=\sum_{c=0}^av^{2c}$ and
$\on{ch}(\on{Hom}(\CO_\bC(-a),\CO_\bC(-b_1)/\CO_\bC(-b_2)))=
\sum_{c=a-b_2+1}^{a-b_1}v^{2c}$.
The part c) is obvious.
\end{proof}

\subsection{}
Let us denote by
$S\chi_{\widetilde{\ul{d}}}=\Lambda^{-1}\chi_{\widetilde{\ul{d}}}$
(resp. $S\chi_{(\widetilde{\ul{d}},\widetilde{\ul{d}}{}')}=
\Lambda^{-1}\chi_{(\widetilde{\ul{d}},\widetilde{\ul{d}}{}')}$)
the character of $\widetilde{T}\times\BC^*$ in the symmetric
algebra $\on{Sym}^\bullet\CT_{\widetilde{\ul{d}}}\fQ_{\ul{d}}$
(resp.
$\on{Sym}^\bullet\CT_{(\widetilde{\ul{d}},\widetilde{\ul{d}}{}')}
\fE_{\ul{d},i}$). It is the inverse of the character of the
corresponding exterior algebra, thus it lies in the fraction field
$\on{Frac}(K^{\widetilde{T}\times\BC^*}(pt))$.

According to the Bott-Lefschetz fixed point formula, the matrix
coefficient
$\bp_*\bq^*_{[\widetilde{\ul{d}}{}',\widetilde{\ul{d}}]}$ of
$\bp_*\bq^*:\ M_{\ul{d}'}\to M_{\ul{d}}$ with respect to the basis
elements $[\widetilde{\ul{d}}]\in K^{\widetilde{T}\times\BC^*}
(\fQ_{\ul{d}}),\ [\widetilde{\ul{d}}{}']\in
K^{\widetilde{T}\times\BC^*} (\fQ_{\ul{d}'})$ (see ~\ref{matrix})
equals $S\chi_{(\widetilde{\ul{d}},\widetilde{\ul{d}}{}')}/
S\chi_{\widetilde{\ul{d}}{}'}$. Similarly, the matrix coefficient
$\bq_*(\fL_i\otimes\bp^*)_{[\widetilde{\ul{d}},\widetilde{\ul{d}}{}']}$
of $\bq_*(\fL_i\otimes\bp^*):\ M_{\ul{d}}\to M_{\ul{d}'}$ equals
$\lambda_{(\widetilde{\ul{d}},\widetilde{\ul{d}}{}')}
S\chi_{(\widetilde{\ul{d}},\widetilde{\ul{d}}{}')}/
S\chi_{\widetilde{\ul{d}}}$.

Hence, the matrix coefficient
$E_{i[\widetilde{\ul{d}},\widetilde{\ul{d}}{}']}$ of $E_i:\
M_{\ul{d}}\to M_{\ul{d}'}$ equals
$$-t_{i+1}^{-i-1}t_i^{i-1}v^{(i-1)d_{i-1}+(i+1)d_{i+1}-2id_i-i}
\lambda_{(\widetilde{\ul{d}},\widetilde{\ul{d}}{}')}
S\chi_{(\widetilde{\ul{d}},\widetilde{\ul{d}}{}')}/
S\chi_{\widetilde{\ul{d}}}.$$ And the matrix coefficient
$F_{i[\widetilde{\ul{d}},\widetilde{\ul{d}}{}']}$ of $F_i:\
M_{\ul{d}}\to M_{\ul{d}'}$ equals
$t_{i+1}^it_i^{-i}v^{2id_i-id_{i-1}-id_{i+1}-i}
S\chi_{(\widetilde{\ul{d}}{}',\widetilde{\ul{d}})}/
S\chi_{\widetilde{\ul{d}}}$.

Thus, Proposition ~\ref{zanudstvo} admits the following Corollary.

\begin{cor}
\label{coefficients}
$$E_{i[\widetilde{\ul{d}},\widetilde{\ul{d}}{}']}=
-t_{i+1}^{-i-1}t_i^{i-1}v^{(i-1)d_{i-1}+(i+1)d_{i+1}-2id_i-i}
t_j^2v^{-2d_{i,j}}\times$$
$$(1-v^2)^{-1}\prod_{j\ne k\leq i}(1-t_j^2t_k^{-2}v^{2d_{i,k}-2d_{i,j}})^{-1}
\prod_{k\leq i-1}(1-t_j^2t_k^{-2}v^{2d_{i-1,k}-2d_{i,j}})$$ if
$d'_{i,j}=d_{i,j}+1$ for certain $j\leq i$;

$$F_{i[\widetilde{\ul{d}},\widetilde{\ul{d}}{}']}=
t_{i+1}^it_i^{-i}v^{2id_i-id_{i-1}-id_{i+1}-i}\times$$
$$(1-v^2)^{-1}\prod_{j\ne k\leq i}(1-t_k^2t_j^{-2}v^{2d_{i,j}-2d_{i,k}})^{-1}
\prod_{k\leq i+1}(1-t_k^2t_j^{-2}v^{2d_{i,j}-2d_{i+1,k}})$$ if
$d'_{i,j}=d_{i,j}-1$ for certain $j\leq i$;

All the other matrix coefficients of $E_i,F_i$ vanish.
\end{cor}

Now the relation ~(\ref{ochevidno}) boils down to the following
identity.

\begin{prop}
\label{mrak}
$$\frac{t_it_{i+1}^{-1}v^{d_{i-1}-2d_i+d_{i+1}-1}-
t_i^{-1}t_{i+1}v^{-d_{i-1}+2d_i-d_{i+1}+1}}{v-v^{-1}}
(1-v^2)^2v^{d_{i-1}-d_{i+1}}t_it_{i+1}=$$
$$\sum_{j\leq i}t_j^2v^{-2d_{i,j}+2}
(1-t_i^2t_j^{-2}v^{2d_{i,j}-2d_{i+1,i}})
(1-t_{i+1}^2t_j^{-2}v^{2d_{i,j}-2d_{i+1,i+1}})\times$$
$$\times\prod_{k\leq i}^{k\ne j}(1-t_k^2t_j^{-2}v^{2d_{i,j}-2d_{i,k}})^{-1}
(1-t_j^2t_k^{-2}v^{2d_{i,k}-2d_{i,j}+2})^{-1}\times$$
$$\times\prod_{k\leq i-1}(1-t_k^2t_j^{-2}v^{2d_{i,j}-2d_{i+1,k}})
(1-t_j^2t_k^{-2}v^{2d_{i-1,k}-2d_{i,j}+2})-$$
$$-\sum_{j\leq i}t_j^2v^{-2d_{i,j}}
(1-t_i^2t_j^{-2}v^{2d_{i,j}-2d_{i+1,i}+2})
(1-t_{i+1}^2t_j^{-2}v^{2d_{i,j}-2d_{i+1,i+1}+2})\times$$
$$\times\prod_{k\leq i}^{k\ne j}(1-t_k^2t_j^{-2}v^{2d_{i,j}-2d_{i,k}+2})^{-1}
(1-t_j^2t_k^{-2}v^{2d_{i,k}-2d_{i,j}})^{-1}\times$$
$$\times\prod_{k\leq i-1}(1-t_k^2t_j^{-2}v^{2d_{i,j}-2d_{i+1,k}+2})
(1-t_j^2t_k^{-2}v^{2d_{i-1,k}-2d_{i,j}}).$$
\end{prop}

\begin{proof} We introduce the new variables $q:=v^2;\ s_j:=t_j^2v^{-2d_{ij}},\
1\leq j\leq i;\ r_k:=t_k^2v^{-2d_{i+1,k}},\ 1\leq k\leq i+1;\
p_k:=t_k^2v^{-2d_{i-1,k}},\ 1\leq k\leq i-1$. Then the LHS of
~\ref{mrak} equals
$$(1-q)\left(q\prod_{k=1}^{i+1}r_k\prod_{j=1}^is_j^{-1}-
\prod_{j=1}^is_j\prod_{k=1}^{i-1}p_k^{-1}\right)$$ while the RHS
of ~\ref{mrak} equals
$$\prod_{j=1}^is_j\prod_{k=1}^{i-1}p_k^{-1}
\left(q\sum_{j\leq i}s_j^{-2}\prod_{k=1}^{i+1}(s_j-r_k)
\prod_{k=1}^{i-1}(p_k-qs_j) \prod_{k\leq i}^{k\ne
j}(s_j-s_k)^{-1}(s_k-qs_j)^{-1}\right.-$$
$$\left.\sum_{j\leq i}s_j^{-2}\prod_{k=1}^{i+1}(s_j-qr_k)
\prod_{k=1}^{i-1}(p_k-s_j) \prod_{k\leq i}^{k\ne
j}(s_j-qs_k)^{-1}(s_k-s_j)^{-1}\right)$$ Dividing both the LHS and
the RHS by $\prod_{j=1}^is_j\prod_{k=1}^{i-1}p_k^{-1}$ we arrive
at
$$(1-q)(q\prod_{j=1}^is_j^{-2}\prod_{k=1}^{i-1}p_k\prod_{k=1}^{i+1}r_k-1)=$$
$$q\sum_{j\leq i}s_j^{-2}\prod_{k=1}^{i+1}(s_j-r_k)
\prod_{k=1}^{i-1}(p_k-qs_j) \prod_{k\leq i}^{k\ne
j}(s_j-s_k)^{-1}(s_k-qs_j)^{-1}-$$
$$\sum_{j\leq i}s_j^{-2}\prod_{k=1}^{i+1}(s_j-qr_k)
\prod_{k=1}^{i-1}(p_k-s_j) \prod_{k\leq i}^{k\ne
j}(s_j-qs_k)^{-1}(s_k-s_j)^{-1}.$$ If we subtract the LHS from the
RHS we obtain a rational expression in $s_j$ of degree 0, that is,
the degree of numerator is not bigger than the degree of
denominator. We see easily that as $s_j$ tends to $\infty$, the
difference of the RHS and the LHS tends to 0. The possible poles
of the difference can occur at $s_j=0,\ s_j=s_k,\ s_j=qs_k,\
s_j=q^{-1}s_k$. We see easily that the principal parts of the
difference at these points vanish. We conclude that the difference
is identically 0. This completes the proof of the Proposition.
\end{proof}

\subsection{}
\label{import} To finish the proof of relation ~(\ref{ochevidno})
we note that the commutator correspondence $E_iF_i-F_iE_i$ is
concentrated on the diagonal of $\fQ_{\ul{d}}\times\fQ_{\ul{d}}$.
This is proved exactly as in Lemma ~\ref{trans}. In other words,
$E_iF_i-F_iE_i$ is given by tensor product
$?\mapsto?\stackrel{L}{\otimes}X_i$ for certain $X_i\in
M_{\ul{d}}$. This means that in the basis $[\widetilde{\ul{d}}]$
the operator $E_iF_i-F_iE_i$ is diagonal. Now the Proposition
~\ref{mrak} computes the matrix coefficient
$(E_iF_i-F_iE_i)_{[\widetilde{\ul{d}},\widetilde{\ul{d}}]}$ and
proves that it equals
$\frac{K_i-K_i^{-1}}{v-v^{-1}}|_{M_{\ul{d}}}$. This completes the
proof of the relation ~(\ref{ochevidno}).

\subsection{}
\label{alternate} Alternatively, the relation ~(\ref{ochevidno})
follows from the next Conjecture. We consider a 2-dimensional
vector space with a basis $\fw_1,\fw_2$. Let $\fT$ be a torus
acting on $\fw_1$ (resp. $\fw_2$) via a character $\tau_1^2$
(resp. $\tau_2^2$). Let $\fZ_{\fd_1,\fd_2}$ be the moduli stack of
flags of coherent sheaves $\fW_1\subset\fW_2$ on $\bC$ locally
free at $\infty\in\bC$, equipped with a trivialization
$\fW_1|_{\infty}=\langle\fw_1\rangle,\
\fW_2|_{\infty}=\langle\fw_1,\fw_2\rangle$, and such that
$\deg\fW_1=-\fd_1,\ \deg\fW_2/\fW_1=-\fd_2$. We have a natural
correspondence
$\fE_{\fd_1}\subset\fZ_{\fd_1,\fd_2}\times\fZ_{\fd_1+1,\fd_2-1}$
formed by the pairs $(\fW_1,\fW_2;\fW'_1,\fW'_2)$ such that
$\fW'_1\subset\fW_1\subset\fW_2=\fW'_2$. The projection
$\fE_{\fd_1}\to\fZ_{\fd_1,\fd_2}$ (resp.
$\fE_{\fd_1}\to\fZ_{\fd_1+1,\fd_2-1}$) is denoted by $\bp$ (resp.
$\bq$). Finally, $\fE_{\fd_1}$ is equipped with the line bundle
$\fL_{\fd_1}$ whose fiber at the point
$(\fW_1,\fW_2;\fW'_1,\fW'_2)$ equals $\Gamma(\bC,\fW_1/\fW'_1)$.

The stack $\fZ_{\fd_1,\fd_2}$ is smooth, and acted upon by
$\fT\times\BC^*$. So it makes sense to consider the operators
$$f:=\bp_*\bq^*:$$ $$K^{\fT\times\BC^*}(\fZ_{\fd_1,\fd_2})
\otimes_{K^{\fT\times\BC^*}(pt)}\on{Frac}(K^{\fT\times\BC^*}(pt))\to
K^{\fT\times\BC^*}(\fZ_{\fd_1-1,\fd_2+1})
\otimes_{K^{\fT\times\BC^*}(pt)}\on{Frac}(K^{\fT\times\BC^*}(pt)),$$
$$e:=-\tau_1^{-1}\tau_2^{-1}v^{\fd_1+\fd_2}
\bq_*(\fL_{\fd_1}\otimes\bp^*):$$
$$K^{\fT\times\BC^*}(\fZ_{\fd_1,\fd_2})
\otimes_{K^{\fT\times\BC^*}(pt)}\on{Frac}(K^{\fT\times\BC^*}(pt))\to
K^{\fT\times\BC^*}(\fZ_{\fd_1+1,\fd_2-1})
\otimes_{K^{\fT\times\BC^*}(pt)}\on{Frac}(K^{\fT\times\BC^*}(pt)),$$
$$K=\tau_1^{-1}\tau_2v^{\fd_1-\fd_2+1}:$$
$$K^{\fT\times\BC^*}(\fZ_{\fd_1,\fd_2})
\otimes_{K^{\fT\times\BC^*}(pt)}\on{Frac}(K^{\fT\times\BC^*}(pt))\to
K^{\fT\times\BC^*}(\fZ_{\fd_1,\fd_2})
\otimes_{K^{\fT\times\BC^*}(pt)}\on{Frac}(K^{\fT\times\BC^*}(pt))$$

\begin{conj}
\label{inteligent} $ef-fe=\frac{K-K^{-1}}{v-v^{-1}}$.
\end{conj}

\subsection{}
\label{derivation} To derive the relation ~(\ref{ochevidno}), or
equivalently, ~(\ref{ochevidno'}) for $j=i$ from Conjecture
~\ref{inteligent} we consider the map
$$\fz_{\ul{d}}:\ \fQ_{\ul{d}}\to\fZ_{d_i-d_{i-1},d_{i+1}-d_i},\
\CW_\bullet\mapsto(\CW_i/\CW_{i-1},\CW_{i+1}/\CW_{i-1}).$$ Then we
have
$$(\fQ_{\ul{d}}\times\fZ_{d_i-d_{i-1}+1,d_{i+1}-d_i-1})
\times_{\fZ_{d_i-d_{i-1},d_{i+1}-d_i}\times\fZ_{d_i-d_{i-1}+1,d_{i+1}-d_i-1}}
\fE_{d_i-d_{i-1}}=\fE_{\ul{d},i}\subset\fQ_{\ul{d}}\times\fQ_{\ul{d}+i}.$$
We also have the natural maps
$$\fe_{\ul{d},i}:\ \fE_{\ul{d},i}\to\fE_{d_i-d_{i-1}},$$
$$^\sT\fe_{\ul{d},i}:\ {}^\sT\fE_{\ul{d},i}\to\ {}^\sT\fE_{d_i-d_{i-1}},$$
$$\fh_{\ul{d},i}:\ \fE_{\ul{d}-i,i}\circ\ {}^\sT\fE_{\ul{d}-i,i}\to
\fE_{d_i-d_{i-1}-1}\circ\ {}^\sT\fE_{d_i-d_{i-1}-1},$$
$$'\fh_{\ul{d},i}:\ {}^\sT\fE_{\ul{d},i}\circ\fE_{\ul{d},i}\to\
{}^\sT\fE_{d_i-d_{i-1}}\circ\fE_{d_i-d_{i-1}}.$$ We may consider
$e_i$ (resp. $f_i,e,f$) as an element of
$K^{\widetilde{T}\times\BC^*}(\fE_{\ul{d},i})$ (resp.
$K^{\widetilde{T}\times\BC^*}(\ {}^\sT\fE_{\ul{d},i}),\\
K^{\widetilde{T}\times\BC^*}(\fE_{d_i-d_{i-1}}),\
K^{\widetilde{T}\times\BC^*}(\ {}^\sT\fE_{d_i-d_{i-1}})$). We
evidently have
$$\fe_{\ul{d},i}^*e=e_i,\ {}^\sT\fe_{\ul{d},i}^*f=f_i.$$
Moreover, according to ~\cite{n}, ~8.2 (Restriction of the
convolution to submanifolds), we have

\begin{equation}
\label{nak} \fh_{\ul{d},i}^*(e*f)=e_i*f_i,\
'\fh_{\ul{d},i}^*(f*e)=f_i*e_i.
\end{equation}

We already know from the argument in ~\ref{import} that the
correspondence $e_i*f_i-f_i*e_i$ acts as tensor multiplication
with a certain class $X_i\in M_{\ul{d}}$. Similarly, the
correspondence $e*f-f*e$ acts in
$K^{\widetilde{T}\times\BC^*}(\fZ_{d_i-d_{i-1},d_{i+1}-d_i})
\otimes_{K^{\widetilde{T}\times\BC^*}(pt)}
\on{Frac}(K^{\widetilde{T}\times\BC^*}(pt))$ as tensor
multiplication with a certain class $\fX\in
K^{\widetilde{T}\times\BC^*}(\fZ_{d_i-d_{i-1},d_{i+1}-d_i})
\otimes_{K^{\widetilde{T}\times\BC^*}(pt)}
\on{Frac}(K^{\widetilde{T}\times\BC^*}(pt))$ By ~(\ref{nak}) we
must have $X_i=\fz_{\ul{d}}^*\fX$. Thus the relation
~(\ref{ochevidno'}) for $j=i$ follows from Conjecture
~\ref{inteligent}.

\subsection{}
\label{univerma} To complete the proof of Theorem ~\ref{main} it
remains to check the relations ~(\ref{Serre1}), ~(\ref{Serre2}).
To this end we consider the algebra $\widetilde{U}$ given by the
generators $E_i,L_i^{\pm1},K_i^{\pm1},F_i,\ 1\leq i\leq n-1$, and
the relations ~(\ref{och})--~(\ref{ochevidno}). Thus, $U$ is the
quotient of $\widetilde{U}$ by the Serre relations.

We extend the scalars to
$\on{Frac}(\BC[L_1^{\pm1},\ldots,L_{n-1}^{\pm1}])$: we set
$$U'=U\otimes_\BC
\on{Frac}(\BC[L_1^{\pm1},\ldots,L_{n-1}^{\pm1}]),\
\widetilde{U}{}'=\widetilde{U}\otimes_\BC
\on{Frac}(\BC[L_1^{\pm1},\ldots,L_{n-1}^{\pm1}])$$ Note that
$\widetilde{U}{}'$ acts in $M$, so $U'$ acts in the quotient
$\overline{M}$ of $M$ by the two-sided ideal $\CI$ in
$\widetilde{U}{}'$ generated by the Serre relations. So it
suffices to check that $\overline{M}=M$, or equivalently, $\CI
M=0$.

Now $M$ has the size of the universal Verma module over $U'$ which
is an irreducible $U'$- (and $\widetilde{U}{}'$-) module. In
effect, a bijection between the set $\{[\widetilde{\ul{d}}]\}$,
and the set of Kostant partitions for $\mathfrak{sl}_n$ is defined
e.g. in ~\cite{fk}, 2.1.1. Hence we only have to check that $\CI
M\ne M$. But any element $x\in\CI$ of principal grading degree 0
annihilates the lowest weight vector $[(0,\ldots,0)]$ of $M$ since
we may shift the generators $e_i$ in the expression of $x$ to the
right.

This completes the proof of the Serre relations in $M$ along with
the proof of Theorem ~\ref{main}.

\subsection{Remark}
\label{joseph} (A.~Joseph) We have constructed a basis
$\{[\widetilde{\ul{d}}]\}$ in the universal Verma module $M$ over
$U$. Though we can not identify it with any known type of basis,
the parametrization of this basis coincides with the polyhedral
realization of the crystal base of $U^+_v(\mathfrak{sl}_n)$
corresponding to the reduced expression in the Weyl group of
$SL_n$:
$$w_0=s_{n-1}s_{n-2}\ldots s_1s_{n-1}s_{n-2}\ldots s_2\ldots s_{n-1}s_{n-2}
s_{n-1}$$ (see ~\cite{nz}).

\subsection{}
\label{Shapovalov} Recall that the universal Verma module $M$ over
$U$ is equipped with the symmetric Shapovalov form $(,)$ with
values in $\on{Frac}(\BC[\widetilde{T}\times\BC^*])$. It is
characterized by the properties

(a) $([\widetilde{\ul{d}}{}_0],[\widetilde{\ul{d}}{}_0])=1$ where
$[\widetilde{\ul{d}}{}_0]=[(0,\ldots,0)]$ is the lowest weight
vector;

(b) $(E_ix,y)=(x,F_iy)\ \ \forall\ x,y\in M$.

\medskip

We will write down a geometric expression for the Shapovalov form.
Evidently, the different weight spaces of $M$ are orthogonal with
respect to the Shapovalov form. We consider the line bundle $\D_i$
on $\fQ_{\ul{d}}$ whose fiber at the point $(\CW_\bullet)$ equals
$\det R\Gamma(\bC,\CW_i)$. We also define the line bundle
$\D_{\ul{d}}:=\bigotimes_{i=1}^{n-1}\D_i$.

\begin{prop}
\label{Shapoval} For $\CG_1,\CG_2\in M_{\ul{d}}$ we have
$$(\CG_1,\CG_2)=(-1)^{\sum_{i=1}^{n-1}d_i}
v^{\sum_{i=1}^{n-1}2id_i^2-\sum_{i=2}^{n-1}(2i-1)d_id_{i-1}}
\prod_{i=1}^nt_i^{(2i-1)(d_{i-1}-d_i)}[R\Gamma(\fQ_{\ul{d}},
\CG_1\otimes\CG_2\otimes\D_{\ul{d}})]$$
\end{prop}

\begin{proof}
Since $\det R\Gamma$ is multiplicative in short exact sequences,
we have an equality of line bundles on the correspondence
$\fE_{\ul{d},i}:\ \bp^*\D_{\ul{d}}=\bq^*\D_{\ul{d}}\otimes\fL_i$. Now the projection
formula shows that the operators $\bp_*\bq^*$ and
$\bq_*(\fL_i\otimes\bp^*)$ are adjoint with respect to the pairing
$\CG_1,\CG_2\mapsto
R\Gamma(\fQ_{\ul{d}},\CG_1\otimes\CG_2\otimes\D_{\ul{d}})$. Finally, it is
easy to see that the $v,t$-factor takes care of the scaling
coefficients of our $E_i,F_i$.
\end{proof}

\subsection{}
\label{conjugate} While the operators $E_i,F_i$ are conjugate to
each other with respect to the Shapovalov form, the operators
$e_i,f_i$ are not. In fact, obviously, $e_i^*=K_i^{2i}f_i$. It is
known that a completion of the universal Verma module $M$ contains
a unique vector $\fk=\sum_{\ul{d}}\fk_{\ul{d}}$ (resp.
$\fw=\sum_{\ul{d}}\fw_{\ul{d}}$) such that
$\fk_{(0,\ldots,0)}=\fw_{(0,\ldots,0)}=[(0,\ldots,0)]$, and
$f_i\fk=(1-v^2)^{-1}\fk$ (resp. $e_i^*\fw=(1-v^2)^{-1}\fw$) for
any $i$ (the {\em Whittaker vectors}).

The following proposition gives a geometric construction of the
Whittaker vectors $\fk,\fw\in M$.

\begin{prop}
\label{hlop} a) $\fk_{\ul{d}}=[\CO_{\ul{d}}]$ (the class of the
structure sheaf of $\fQ_{\ul{d}}$);

b) $\fw_{\ul{d}}=
v^{\sum_{i=1}^{n-1}(1-2i)d_i^2-\sum_{i=2}^{n-1}(2-2i)d_id_{i-1}-
\sum_{i=1}^{n-1}d_i}
\prod_{i=1}^nt_i^{(2-2i)(d_{i-1}-d_i)}[\D^{-1}_{\ul{d}}]$.
\end{prop}

\begin{proof}
a) We have $\bq^*\CO_{\ul{d}+i}=\CO_{\fE_{\ul{d},i}}$.
Furthermore, since $\bp\times\br:\
\fE_{\ul{d},i}\to\fQ_{\ul{d}}\times(\bC-\infty)$ is proper and
birational, and both the source and the target are smooth, we have
$(\bp\times\br)_*[\CO_{\fE_{\ul{d},i}}]=
[\CO_{\ul{d}}]\boxtimes[\CO_{\bC-\infty}]$. In effect,
$(\bp\times\br)_*\CO_{\fE_{\ul{d},i}}=
\CO_{\ul{d}}\boxtimes\CO_{\bC-\infty}$, and the higher direct
images $R^{>0}(\bp\times\br)_*\CO_{\fE_{\ul{d},i}}$ vanish.
Finally, $pr_*[\CO_{\ul{d}}\boxtimes\CO_{\bC-\infty}]=
(1-v^2)^{-1}[\CO_{\ul{d}}]$ where $pr:\
\fQ_{\ul{d}}\times(\bC-\infty)\to \fQ_{\ul{d}}$ is the projection
to the first factor.

b) Recall that $e_i^*=K_i^{2i}f_i$. Thus we have to check that
$f_i[\D^{-1}_{\ul{d}+i}]=t_i^2v^{2d_{i-1}-2d_i}(1-v^2)^{-1}[\D^{-1}_{\ul{d}}]$.
Furthermore, recall that on $\fE_{\ul{d},i}$ we have a canonical
isomorphism
$\bq^*\D^{-1}_{\ul{d}+i}=\fL_i\otimes\bp^*\D^{-1}_{\ul{d}}$. By
the projection formula we are reduced to
\begin{equation}
\label{zhelaem}
\bp_*[\fL_i]=t_i^2v^{2d_{i-1}-2d_i}(1-v^2)^{-1}[\CO_{\ul{d}}]
\end{equation}
This can be calculated in the basis $[\widetilde{\ul{d}}]$ where
we already know the matrix coefficients of our operators (see
Corollary ~\ref{coefficients}). More precisely, by the
Bott-Lefschetz fixed point formula, we have to check
$$\sum_{j\leq i}t_j^2v^{-2d_{ij}}(1-v^2)^{-1}
\prod_{j\ne k\leq i}(1-t_j^2t_k^{-2}v^{2d_{i,k}-2d_{i,j}})^{-1}
\prod_{k\leq i-1}(1-t_j^2t_k^{-2}v^{2d_{i-1,k}-2d_{i,j}})=$$
$$=t_i^2v^{2d_{i-1}-2d_i}(1-v^2)^{-1}$$
Recall the change of variables we used in the proof of Proposition
~\ref{mrak}: $s_j:=t_j^2v^{-2d_{ij}},\ 1\leq j\leq i;\
p_k:=t_k^2v^{-2d_{i-1,k}},\ 1\leq k\leq i-1$. Then we have to
prove
$$\sum_{j\leq i}s_j\prod_{k\leq i}^{k\ne j}(1-s_js_k^{-1})^{-1}
\prod_{k\leq i-1}(1-s_jp_k^{-1})=s_1\cdots s_ip_1^{-1}\cdots
p_{i-1}^{-1}$$ This follows immediately from the well known
identity
$$\sum_{j\leq i}\prod_{k\leq i-1}(p_k-s_j)
\prod_{k\leq i}^{k\ne j}(s_k-s_j)^{-1}=1.$$ This completes the
proof of the Proposition.
\end{proof}

\begin{cor}
\label{vot tebe} The Shapovalov scalar product of the Whittaker
vectors equals
$(\fk_{\ul{d}},\fw_{\ul{d}})=(-1)^{\sum_{i=1}^{n-1}d_i}
v^{\sum_{i=1}^{n-1}d_i^2-\sum_{i=2}^{n-1}d_id_{i-1}-\sum_{i=1}^{n-1}d_i}
\prod_{i=1}^nt_i^{d_{i-1}-d_i}[R\Gamma(\fQ_{\ul{d}},\CO_{\ul{d}})]$.
\end{cor}

\subsection{}
\label{generating} According to the works ~\cite{e}, ~\cite{s2},
the appropriate generating function of the Shapovalov scalar
product of the Whittaker vectors satisfies a $v$-deformed
($v$-difference) version of the quantum Toda lattice equations.
Let us recall the required notations and results.

We introduce the formal variables $\sz_1,\ldots,\sz_n$, and we set
$\sQ_i=\exp(\sz_i-\sz_{i+1}),\ i=1,\ldots,n-1$. We set
$\hbar=\log(v)$, so that $v=\exp(\hbar)$. We introduce the shift
operators $\sT_i,\ i=1,\ldots,n$, acting on the space of functions
of $\sz_1,\ldots,\sz_n$ invariant with respect to the simultaneous
translations
$f(\sz_1,\ldots,\sz_n)=f(\sz_1+\sz,\ldots,\sz_n+\sz)$. Namely, we
set $\sT_if(\sz_1,\ldots,\sz_n)=
f(\sz_1,\ldots,\sz_i+\hbar,\ldots,\sz_n)$.

We define the following $v$-difference operators:

\begin{equation}
\label{Etingof} {\mathfrak
S}:=\sum_{j=1}^n\sT_j^2+v^{-2}\sum_{i=1}^{n-1}\sQ_i\sT_i\sT_{i+1}
\end{equation}

\begin{equation}
\label{Givental} {\mathfrak
G}:=\sT_1^2+\sT_2^2(1-\sQ_1)+\ldots+\sT_n^2(1-\sQ_{n-1})
\end{equation}

We also consider the following generating functions:

\begin{equation}
\label{etingof} {\mathfrak
I}:=\prod_{i=1}^{n-1}\sQ_i^{\frac{-\log(t_1\cdots t_i)}{\hbar}}
\sum_{\ul{d}}(\fk_{\ul{d}},\fw_{\ul{d}})\sQ_1^{d_1}\cdots\sQ_{n-1}^{d_{n-1}}
\end{equation}

\begin{equation}
\label{givental} {\mathfrak
J}:=\prod_{i=1}^{n-1}\sQ_i^{\frac{-\log(t_1\cdots t_i)}{\hbar}}
\sum_{\ul{d}}[R\Gamma(\fQ_{\ul{d}},
\CO_{\ul{d}})]\sQ_1^{d_1}\cdots\sQ_{n-1}^{d_{n-1}}
\end{equation}

Then according to the last formula of ~\cite{s2} (or equivalently,
the formula ~(5.7) of ~\cite{e}), we have

\begin{equation}
\label{etisev} {\mathfrak S}{\mathfrak
I}=\left(\sum_{i=1}^nt_i^2\right){\mathfrak I}
\end{equation}
In effect, the seeming discrepancy between the formula
~(\ref{Etingof}) above, and the formula ~(5.7) of ~\cite{e} is
explained by the fact that (a) our $v$ corresponds to $q$ of
~\cite{e}; (b) our Whittaker vectors have eigenvalue
$(1-v^2)^{-1}$, whereas the Whittaker vectors of ~\cite{e} have
eigenvalue 1, which takes care of the factor $(q-q^{-1})^2$ in the
second summand of the formula ~(5.7) of ~\cite{e}.

Now the argument of ~\cite{e}, section ~6 (see the formula ~(6.5))
together with Corollary ~\ref{vot tebe}, establishes

\begin{equation}
\label{talgiven} {\mathfrak G}{\mathfrak
J}=\left(\sum_{i=1}^nt_i^2\right){\mathfrak J}
\end{equation}

thus reproving the Main Theorem ~2 of ~\cite{gl}.

\section{Parabolic sheaves and affine quantum groups}
In this section we want to generalize the previous results to the
affine setting. \label{p}

\subsection{Parabolic sheaves}
We recall the setup of ~\cite{fgk}. Let $\bX$ be another smooth
projective curve of genus zero. We fix a coordinate $x$ on $\bX$,
and consider the action of $\BC^*$ on $\bX$ such that
$u(x)=u^{-2}x$. We have $\bX^{\BC^*}=\{0_\bX,\infty_\bX\}$. Let
$\bS$ denote the product surface $\bC\times\bX$. Let $\bD_\infty$
denote the divisor $\bC\times\infty_\bX\cup\infty_\bC\times\bX$.
Let $\bD_0$ denote the divisor $\bC\times0_\bX$.

Given an $n$-tuple of nonnegative integers
$\ul{d}=(d_0,\ldots,d_{n-1})$, we say that a {\em parabolic sheaf}
$\CF_\bullet$ of degree $\ul{d}$ is an infinite flag of torsion
free coherent sheaves of rank $n$ on $\bS:\
\ldots\subset\CF_{-1}\subset\CF_0\subset\CF_1\subset\ldots$ such
that:

(a) $\CF_{k+n}=\CF_k(\bD_0)$ for any $k$;

(b) $ch_1(\CF_k)=k[\bD_0]$ for any $k$: the first Chern classes
are proportional to the fundamental class of $\bD_0$;

(c) $ch_2(\CF_k)=d_i$ for $i\equiv k\pmod{n}$;

(d) $\CF_0$ is locally free at $\bD_\infty$ and trivialized at
$\bD_\infty:\ \CF_0|_{\bD_\infty}=W\otimes\CO_{\bD_\infty}$;

(e) For $-n\leq k\leq0$ the sheaf $\CF_k$ is locally free at
$\bD_\infty$, and the quotient sheaves $\CF_k/\CF_{-n},\
\CF_0/\CF_k$ (both supported at $\bD_0=\bC\times0_\bX\subset\bS$)
are both locally free at the point $\infty_\bC\times0_\bX$;
moreover, the local sections of $\CF_k|_{\infty_\bC\times \bX}$
are those sections of $\CF_0|_{\infty_\bC\times
\bX}=W\otimes\CO_\bX$ which take value in $\langle
w_1,\ldots,w_{n-k}\rangle\subset W$ at $\infty_\bX\in \bX$.

\medskip

According to ~\cite{fgk}, ~3.5, the fine moduli space
$\CP_{\ul{d}}$ of degree $\ul{d}$ parabolic sheaves exists and is
a smooth connected quasiprojective variety of dimension
$2d_0+\ldots+2d_{n-1}$.

The group $\widetilde{T}\times\BC^*\times\BC^*$ acts naturally on
$\CP_{\ul{d}}$, and its fixed point set is finite.

\subsection{Correspondences}
If the collections $\ul{d}$ and $\ul{d}'$ differ at the only place
$i\in I:=\BZ/n\BZ$, and $d'_i=d_i+1$, then we consider a
correspondence
$\sE_{\ul{d},i}\subset\CP_{\ul{d}}\times\CP_{\ul{d}'}$ formed by
the pairs $(\CF_\bullet,\CF'_\bullet)$ such that for $j\not\equiv
i\pmod{n}$ we have $\CF_j=\CF'_j$, and for $j\equiv i\pmod{n}$ we
have $\CF'_j\subset\CF_j$.

It is a smooth quasiprojective algebraic variety of dimension
$2\sum_{i\in I}d_i+1$. In effect, the argument of ~\cite{fgk},
~Lemma~3.3, reduces this statement to the corresponding fact about
Laumon correspondences (see ~\cite{la1}, ~2.10).

We denote by $\bp$ (resp. $\bq$) the natural projection
$\sE_{\ul{d},i}\to\CP_{\ul{d}}$ (resp.
$\sE_{\ul{d},i}\to\CP_{\ul{d}'}$). For $j\equiv i\pmod{n}$ the
correspondence $\sE_{\ul{d},i}$ is equipped with a natural line
bundle $\sL_j$ whose fiber at $(\CF_\bullet,\CF'_\bullet)$ equals
$\Gamma(\bC,\CF_{j-n}/\CF'_{j-n})$. Finally, we have a transposed
correspondence
$^\sT\sE_{\ul{d},i}\subset\CP_{\ul{d}'}\times\CP_{\ul{d}}$.

\subsection{}
We denote by ${}'\CM$ the direct sum of equivariant (complexified)
$K$-groups:
${}'\CM=\oplus_{\ul{d}}K^{\widetilde{T}\times\BC^*\times\BC^*}(\CP_{\ul{d}})$.
It is a module over $K^{\widetilde{T}\times\BC^*\times\BC^*}(pt)
=\BC[\widetilde{T}\times\BC^*\times\BC^*]= \BC[t_1,\ldots,t_n,v,u\
:\ t_1\cdots t_n=1]$. We define $\CM=\
{}'\CM\otimes_{K^{\widetilde{T}\times\BC^*\times\BC^*}(pt)}
\on{Frac}(K^{\widetilde{T}\times\BC^*\times\BC^*}(pt))$.

We have an evident grading $\CM=\oplus_{\ul{d}}\CM_{\ul{d}},\
\CM_{\ul{d}}=K^{\widetilde{T}\times\BC^*\times\BC^*}(\CP_{\ul{d}})
\otimes_{K^{\widetilde{T}\times\BC^*\times\BC^*}(pt)}
\on{Frac}(K^{\widetilde{T}\times\BC^*\times\BC^*}(pt))$.

\subsection{}
\label{operators'} The grading and the correspondences
$^\sT\sE_{\ul{d},i},\sE_{\ul{d},i}$ give rise to the following
operators on $\CM$ (note that though $\bp$ is not proper, $\bp_*$
is well defined on the localized equivariant $K$-theory due to the
finiteness of the fixed point sets):

$K_i=t_{i+1}t_i^{-1}u^{\delta_{0,i}}v^{2d_i-d_{i-1}-d_{i+1}+1}:\
\CM_{\ul{d}}\to\CM_{\ul{d}}$,

$C=uv^n$,

For $i=0,\ldots,n-1$ we define $L_i=t_1^{-1}\cdots
t_i^{-1}v^{d_i+\frac{1}{2}i(n-i)}:\ \CM_{\ul{d}}\to \CM_{\ul{d}}$
(that is, $L_0=v^{d_0}$),

$f_i=\bp_*\bq^*:\ \CM_{\ul{d}}\to\CM_{\ul{d}-i}$;

For $n>2$ and $i=0,\ldots,n-1$ we define
$F_i=t_{i+1}^{-1}v^{d_{i+1}-d_i+\frac{n+1}{2}-i}\bp_*\bq^*:\
\CM_{\ul{d}}\to\CM_{\ul{d}-i}$,

For $n=2$ we define $F_i=f_i$,

$e_i=-t_i^{-1}t_{i+1}^{-1}u^{\delta_{0,i}}v^{d_{i+1}-d_{i-1}}
\bq_*(\sL_i\otimes\bp^*):\ \CM_{\ul{d}}\to\CM_{\ul{d}+i}$,

For $n>2$ and $i=0,\ldots,n-1$ we define
$E_i=-t_i^{-1}u^{\delta_{0,i}}v^{d_i-d_{i-1}+\frac{1-n}{2}+i}
\bq_*(\sL_i\otimes\bp^*):\ \CM_{\ul{d}}\to\CM_{\ul{d}+i}$,

For $n=2$ we define $E_i=e_i$.

\subsection{Sevostyanov's form of affine quantum $SL_n$}
Let $I$ denote the set $\BZ/n\BZ$ of residue classes modulo $n$.

$\CU$ is the $\BC[v,v^{-1}]$-algebra with generators
$E_i,L_i^{\pm1},K_i^{\pm1},C^{\pm1},F_i,\ i\in\BZ/n\BZ$, subject
to the following relations:

\begin{equation}
\label{ev} L_iL_j=L_jL_i,\ K_i=L_i^2L_{i+1}^{-1}L_{i-1}^{-1}C^{\delta_{i,0}}
\end{equation}

\begin{equation}
\label{evid} L_jE_iL_j^{-1} =v^{\delta_{i,j}}E_i,\ L_jF_iL_j^{-1}
=v^{-\delta_{i,j}}F_i,\
C\hphantom{m}\on{is}\hphantom{m}\on{central}
\end{equation}

\begin{equation}
\label{evident}
E_iF_j-F_jE_i=\delta_{i,j}\frac{K_i-K_i^{-1}}{v-v^{-1}}
\end{equation}

\begin{equation}
\label{Ser1} |i-j|>1\ \Longrightarrow\
E_iE_j-E_jE_i=0=F_iF_j-F_jF_i
\end{equation}

\begin{equation}
\label{Ser2} n>2\ \&\ |i-j|=1\ \Longrightarrow\
E_i^2E_j-(v+v^{-1})E_iE_jE_i+E_jE_i^2=0=
F_i^2F_j-(v+v^{-1})F_iF_jF_i+F_jF_i^2
\end{equation}

\begin{equation}
\label{Ser22} n=2\ \&\ |i-j|=1\ \Longrightarrow\
E_i^3E_j-(v^2+1+v^{-2})E_i^2E_jE_i+
(v^2+1+v^{-2})E_iE_jE_i^2-E_jE_i^3=0
\end{equation}

\begin{equation}
\label{Ser22F} n=2\ \&\ |i-j|=1\ \Longrightarrow\
F_i^3F_j-(v^2+1+v^{-2})F_i^2F_jF_i+
(v^2+1+v^{-2})F_iF_jF_i^2-F_jF_i^3=0
\end{equation}

For $n>2$ we also consider elements $e_i,f_i\in\CU$ depending on
the following choice of $n\times n$-matrices $n_{ij},\ c_{ij}$
(cf. ~\cite{s}, ~Remark ~3):

\begin{equation}
\label{choice} n_{i,i}=1,\ n_{i,i+1}=-1,\ n_{i+1,i}=0,
\end{equation}
otherwise $n_{ij}=0$.

\begin{equation}
\label{Russia} c_{i,i+1}=-1,\ c_{i+1,i}=1,
\end{equation}
otherwise $c_{ij}=0$.

Then we set
\begin{equation}
\label{Sevo} f_i:=L_iL_{i+1}^{-1}F_i,\ e_i:=E_iL_i^{-1}L_{i+1}.
\end{equation}

Clearly, the algebra $\CU$ is generated by
$e_i,L_i^{\pm1},K_i^{\pm1}, C^{\pm1},f_i,\ i\in\BZ/n\BZ$, and the
relations ~(\ref{evid})--~(\ref{Ser2}) above are equivalent to the
relations ~(\ref{evid'})--~(\ref{Ser2'}) below.

\begin{equation}
\label{evid'} L_je_iL_j^{-1} =v^{\delta_{i,j}}e_i,\ L_jf_iL_j^{-1}
=v^{-\delta_{i,j}}f_i,\
C\hphantom{m}\on{is}\hphantom{m}\on{central}
\end{equation}

\begin{equation}
\label{evident'}
e_if_j-v^{c_{ij}}f_je_i=\delta_{i,j}\frac{K_i-K_i^{-1}}{v-v^{-1}}
\end{equation}

\begin{equation}
\label{Ser1'} |i-j|>1\ \Longrightarrow\
e_ie_j-e_je_i=0=f_if_j-f_jf_i
\end{equation}

\begin{equation}
\label{Ser2'} |i-j|=1\ \Longrightarrow\
e_i^2e_j-v^{c_{ij}}(v+v^{-1})e_ie_je_i+v^{2c_{ij}}e_je_i^2=0=
f_i^2f_j-v^{c_{ij}}(v+v^{-1})f_if_jf_i+v^{2c_{ij}}f_jf_i^2
\end{equation}

\subsection{} The following
is an affine analogue of Theorem ~\ref{main}. Recall the operators
$E_i,e_i,K_i^{\pm1},L_i^{\pm1},C^{\pm1},F_i,f_i,\ i\in I$, on
$\CM$ defined in ~\ref{operators'}.

\begin{conj}
\label{main'} The operators
$E_i,K_i^{\pm1},L_i^{\pm1},C^{\pm1},F_i,\ i\in I$, on $\CM$
satisfy the relations ~(\ref{ev})--(\ref{Ser22F}). Equivalently,
if $n>2$, the operators $e_i,K_i^{\pm1},L_i^{\pm1},C,f_i,\ i\in
I$, satisfy the relations ~(\ref{ev}),
~(\ref{evid'})--(\ref{Ser2'}).
\end{conj}

\subsection{}
\label{long} We can prove Conjecture ~\ref{main'} for $n>2$ . Let
us sketch  this proof. It is parallel to the proof of Theorem
~\ref{main}. In effect, the relation ~(\ref{evident'}) for $i\ne
j$ follows from the transversality statement absolutely similar to
Lemma ~\ref{trans}. More precisely, the argument of ~\cite{fgk}
~(Lemma ~3.3), reduces the required smoothness to that proved in
Lemma ~\ref{trans}.

The relation ~(\ref{evident'}) for $j=i$ follows from Conjecture
~\ref{inteligent} by the argument of ~\ref{derivation}. Since we
can not prove Conjecture ~\ref{inteligent} at the moment, we will
derive the relation ~(\ref{evident'}) for $j=i$ from its weaker
but accessible form.

To this end we consider the following closed substack
$\fZ'_{\fd_1,\fd_2}\subset\fZ_{\fd_1,\fd_2}$. Recall that a
coherent sheaf $\fW_1$ (resp. $\fW_2$) contains the maximal
torsion subsheaf $\fW_1^{tors}$ (resp. $\fW_2^{tors}$) with the
locally free quotient sheaf $\fW_1^{free}$ (resp. $\fW_2^{free}$).
Moreover, we have $\fW_1\simeq\fW_1^{tors}\oplus\fW_1^{free}$
(resp. $\fW_2\simeq\fW_2^{tors}\oplus\fW_2^{free}$). The closed
substack $\fZ'_{\fd_1,\fd_2}\subset\fZ_{\fd_1,\fd_2}$ classifies
the flags of coherent sheaves (with trivialization at
$\infty\in\bC$) $\fW_1\subset\fW_2$ such that
$\deg\fW_1^{free}\leq0\geq\deg\fW_2^{free}$. We define
$K^{\fT\times\BC^*}(\fZ'_{\fd_1,\fd_2})$ as the $K$-group of
$\fT\times\BC^*$-equivariant coherent sheaves on the smooth stack
$\fZ_{\fd_1,\fd_2}$ supported on the closed substack
$\fZ'_{\fd_1,\fd_2}$. Note that for any
$\ul{d}=(d_1,\ldots,d_{n-1})$ the map $\fz_{\ul{d}}:\
\fQ_{\ul{d}}\to\fZ_{d_i-d_{i-1},d_{i+1}-d_i}$ factors through the
same named map into the closed substack
$\fZ'_{d_i-d_{i-1},d_{i+1}-d_i}$. Similarly, for any
$\ul{d}=(d_0,d_1,\ldots,d_{n-1})$ the map
$$\fz\fz_{\ul{d}}:\ \CP_{\ul{d}}\to\fZ_{d_i-d_{i-1},d_{i+1}-d_i},\
\CF_\bullet\mapsto(\CF_i/\CF_{i-1},\CF_{i+1}/\CF_{i-1})$$ factors
through the same named map into the closed substack
$\fZ'_{d_i-d_{i-1},d_{i+1}-d_i}$.

Let $(\fW_1\subset\fW_2)$ be a $\fT\times\BC^*$-fixed point of
$\fZ'_{\fd_1,\fd_2}$. Let $\iota_{(\fW_1\subset\fW_2)}$ denote its
locally closed embedding into $\fZ'_{\fd_1,\fd_2}$. Let
$Aut_{(\fW_1\subset\fW_2)}$ stand for its automorphisms' group.
One can easily check the following

\begin{lem}
\label{easily} There exists $n,\ i,\ 1\leq i\leq n-1,\
\ul{d}=(d_1,\ldots,d_n)$, and a fixed point
$\widetilde{\ul{d}}\in\fQ_{\ul{d}}^{\widetilde{T}\times\BC^*}$
such that

(a) $\fz_{\ul{d}}(\widetilde{\ul{d}})=(\fW_1\subset\fW_2)$;

(b) $\fz_{\ul{d}}(\widetilde{T}\times\BC^*)$ is a maximal torus of
$Aut_{(\fW_1\subset\fW_2)}$.
\end{lem}

\subsection{}
\label{birka} One way to prove Conjecture ~\ref{inteligent} would
be to reverse the argument of ~\ref{derivation} and derive it from
the relations ~(\ref{ochevidno'}) for all $n,i$. In effect, we
must compute (notations of ~\ref{derivation}) $\fX\in
K^{\fT\times\BC^*}(\fZ_{\fd_1,\fd_2})
\otimes_{K^{\fT\times\BC^*}(pt)}
\on{Frac}(K^{\fT\times\BC^*}(pt))$ while we know
$\fz_{\ul{d}}^*\fX$ for all $n,i,\ul{d}$ such that
$d_i-d_{i-1}=\fd_1,\ d_{i+1}-d_i=\fd_2$ (also, the homomorphism of
tori $\widetilde{T}_n\to\fT$ acts on the characters as
$\tau_1=t_i,\ \tau_2=t_{i+1}$).

Let us denote by $\fY\in K^{\fT\times\BC^*}(\fZ'_{\fd_1,\fd_2})
\otimes_{K^{\fT\times\BC^*}(pt)}
\on{Frac}(K^{\fT\times\BC^*}(pt))$ the restriction of $\fX$ to
$\fZ'_{\fd_1,\fd_2}$.

The Lemma ~\ref{easily} implies that the kernel $Ker_1$ of the
direct product of inverse images
$$\prod_{n,i,\ul{d}}\fz_{\ul{d}}^*:\ K^{\fT\times\BC^*}(\fZ'_{\fd_1,\fd_2})
\otimes_{K^{\fT\times\BC^*}(pt)}\on{Frac}(K^{\fT\times\BC^*}(pt))\to
\prod_{n,i,\ul{d}}M_{\ul{d}}$$ coincides with the kernel $Ker_2$
of the direct product of restrictions
$$\prod_{(\fW_1\subset\fW_2)\in(\fZ'_{\fd_1,\fd_2})^{\fT\times\BC^*}}
\iota^*_{(\fW_1\subset\fW_2)}:\
K^{\fT\times\BC^*}(\fZ'_{\fd_1,\fd_2})
\otimes_{K^{\fT\times\BC^*}(pt)}\on{Frac}(K^{\fT\times\BC^*}(pt))\to$$
$$\to\prod_{(\fW_1\subset\fW_2)\in(\fZ'_{\fd_1,\fd_2})^{\fT\times\BC^*}}
K^{\fT\times\BC^*\times Aut_{(\fW_1\subset\fW_2)}}(pt)
\otimes_{K^{\fT\times\BC^*}(pt)}\on{Frac}(K^{\fT\times\BC^*}(pt))$$
It follows that for any $n,\ i,\ 0\leq i\leq n-1,\
\ul{d}=(d_0,\ldots,d_n)$, such that $\fd_1=d_i-d_{i-1},\
\fd_2=d_{i+1}-d_i$, the kernel $Ker_1=Ker_2$ is contained in the
kernel $Ker_3$ of the inverse image
$$\fz\fz_{\ul{d}}^*:\ K^{\fT\times\BC^*}(\fZ'_{\fd_1,\fd_2})
\otimes_{K^{\fT\times\BC^*}(pt)}\on{Frac}(K^{\fT\times\BC^*}(pt))\to
\CM_{\ul{d}}$$

By the argument of ~\ref{derivation} we know that
$\fY=\frac{K-K^{-1}}{v-v^{-1}}$ modulo $Ker_1$, and hence the same
holds modulo $Ker_3$. The argument of {\em loc. cit.} then shows
that the relation ~(\ref{evident'}) for $j=i$ holds in $\CM$.

\subsection{}
\label{longer} It remains to check the Serre relations. The
relations for negative generators follow from the relations for
positive generators because they are adjoint with respect to the
nondegenerate Shapovalov form, see ~\ref{Shapo} below. So it
suffices to consider the relations ~(\ref{Ser1}), ~(\ref{Ser2})
between $E_i,E_j,\ i\ne j$. It is here that we need the assumption
$n>2$ for technical reasons. Namely, for $n>2$ we can find $k\in
I$ such that $i\ne k\ne j$.

We consider an $n$-dimensional vector space with a basis
$\fw_1,\ldots,\fw_n$, and a torus $\fT$ acting on $\fw_l$ by the
character $\tau^2_l$. Let $\fZ_n$ be the moduli stack of flags of
coherent sheaves $\fW_1\subset\ldots\subset\fW_n$ on $\bC$ locally
free at $\infty\in\bC$, equipped with compatible trivializations
$\fW_l|_\infty=\langle\fw_1,\ldots,\fw_l\rangle$. Note that
$\fZ_n$ has connected components numbered by the degrees of
$\fW_l$, which for $n=2$ coincide with the stacks
$\fZ_{\fd_1,\fd_2}$. Absolutely similarly to ~\ref{alternate} we
introduce the correspondences between various connected
components, which give rise to the operators
$E^\fZ_1,\ldots,E^\fZ_{n-1}$ on the localized equivariant
$K$-theory of $\fZ_n$.

As in ~\ref{long} above, we have a closed substack
$\fZ'_n\subset\fZ_n$ classifying the flags such that
$\deg\fW_l^{free}\leq0,\ 1\leq l\leq n$.

We have a map
$$\fz\fz_k:\ \CP_{\ul{d}}\to\fZ_n,\ (\CF_\bullet)\mapsto
(\CF_{k+1}/\CF_k\subset\ldots\subset\CF_{k+n}/\CF_k)$$ factoring
through the same named map $\CP_{\ul{d}}\to\fZ'_n$. For any $N\geq
n$, and $m$ such that $0\leq m\leq N-n$, and
$\ul{d}=(d_1,\ldots,d_N)$, we also have a map
$$\fz_{m,\ul{d}}:\ \fQ_{\ul{d}}\to\fZ_n,\
(\CW_\bullet)\mapsto(\CW_{m+1}/\CW_m\subset\ldots\subset\CW_{m+n}/\CW_m)$$
factoring through the same named map $\fQ_{\ul{d}}\to\fZ'_n$.

Now the argument of ~\ref{derivation} shows that the Serre
relation between $E_i,E_j$ would follow from the Serre relation
between $E^\fZ_{i'},E^\fZ_{j'}$ for certain $i',j'$. Though we
cannot establish the latter relations, the argument of
~\ref{birka} shows that they hold modulo the subspace $Ker_1$
(because we already know the Serre relations for $\mathfrak{sl}_N$
with arbitrary $N$), and also shows that this suffices to derive
the former relations.

This completes the proof of the Serre relations for $n>2$. Thus,
Conjecture ~\ref{main'} is proved for $n>2$.

\subsection{}
Similarly to ~\ref{Shapoval}, we will write down a geometric
expression for a Shapovalov form on $\CM$, that is a symmetric
$\operatorname{Frac}(\BC[\widetilde{T}\times\BC^*\times\BC^*])$-valued
bilinear form on $\CM$ such that $(E_im_1,m_2)=(m_1,F_im_2)$ for
any $i\in I$, and $m_1,m_2\in\CM$. The different weight spaces of
$\CM$ will be orthogonal with respect to this geometric Shapovalov
form. For $i=0,\ldots,n-1$, we consider the line bundle $\D_i$ on
$\CP_{\ul{d}}$ whose fiber at the point $(\CF_\bullet)$ equals
$\det R\Gamma(\bS,\CF_{i-n})$. We also define the line bundle
$\D_{\ul{d}}:=\bigotimes_{i=0}^{n-1}\D_i$. For
$\CG_1,\CG_2\in\CM_{\ul{d}}$, we set
\begin{equation}
\label{Chapo} (\CG_1,\CG_2):=(-1)^{\sum_{i=0}^{n-1}d_i}
v^{-\sum_{i=0}^{n-1}d_i^2+\sum_{i=0}^{n-1}d_id_{i+1}
+\sum_{i=0}^{n-1}(n-2i)d_i}u^{-d_0}
\prod_{i=0}^{n-1}t_i^{d_i-d_{i-1}}R\Gamma(\CP_{\ul{d}},
\CG_1\otimes\CG_2\otimes\D_{\ul{d}})
\end{equation}
Clearly, the form $(,)$ is nondegenerate, since the classes of the
structure sheaves of the
$\widetilde{T}\times\BC^*\times\BC^*$-fixed points form an
orthogonal basis of $\CM$.

The following proposition is proved exactly as~Proposition~\ref{Shapoval}.

\begin{prop}
\label{Shapo} For $i\in I,\ \CG_1\in\CM_{\ul{d}},\
\CG_2\in\CM_{\ul{d}+i}$ we have
$(E_i\CG_1,\CG_2)=(\CG_1,F_i\CG_2)$.
\end{prop}

\subsection{}
We define a formal sum in a completion of $\CM$ as follows:
$\fn=\sum_{\ul{d}}\fn_{\ul{d}}:=\sum_{\ul{d}}[\CO_{\ul{d}}]
=\sum_{\ul{d}}[\CO_{\CP_{\ul{d}}}]$. We also consider the
following formal sum: $\fu=\sum_{\ul{d}}\fu_{\ul{d}}$ where

\begin{equation}
\fu_{\ul{d}}=v^{2\sum_{i=0}^{n-1}d_i^2-2\sum_{i=0}^{n-1}d_id_{i+1}
-\sum_{i=0}^{n-1}(n-2i+1)d_i}u^{2d_0}
\prod_{i=1}^nt_i^{2d_{i-1}-2d_i}[\D^{-1}_{\ul{d}}]
\end{equation}

\begin{prop}
\label{vykusi} a) $\fn$ is a common eigenvector of the operators
$f_i$ with the eigenvalue $(1-v^2)^{-1}$;

b) $\fu$ is a common eigenvector of the operators $e_i^*$ with the
eigenvalue $(1-v^2)^{-1}$.
\end{prop}

\begin{proof} a) is proved exactly as Proposition ~\ref{hlop} ~(a).

To check b) we argue as in the proof of Proposition ~\ref{hlop}
~(b), and reduce it to
\begin{equation}
\label{hotim}
\bp_*[\sL_i]=t_i^2u^{-2\delta_{0,i}}v^{2d_{i-1}-2d_i}
(1-v^2)^{-1}[\CO_{\ul{d}}]
\end{equation}
To verify this we recall the setup of ~\ref{alternate}, and claim
that in the notations of {\em loc. cit.} we have
\begin{equation}
\label{imeem} \bp_*[\fL_{\fd_1}]=\tau_1^2v^{-2\fd_1}(1-v^2)^{-1}
[\CO_{\fZ_{\fd_1,\fd_2}}]
\end{equation}
In effect, ~(\ref{imeem}) is deduced from ~(\ref{zhelaem}) by the
argument of ~\ref{birka}. Finally, ~(\ref{hotim}) is deduced from
~(\ref{imeem}) by the argument of ~\ref{derivation}.

The Proposition is proved.
\end{proof}

\begin{cor}
\label{vot te} The Shapovalov scalar product of the Whittaker
vectors equals
$(\fn_{\ul{d}},\fu_{\ul{d}})=(-1)^{\sum_{i=0}^{n-1}d_i}
v^{\sum_{i=0}^{n-1}d_i^2-\sum_{i=0}^{n-1}d_id_{i-1}-\sum_{i=0}^{n-1}d_i}
u^{d_0}\prod_{i=1}^nt_i^{d_{i-1}-d_i}R\Gamma(\CP_{\ul{d}},\CO_{\ul{d}})$.
\end{cor}

\subsection{}
\label{new} We define $\CM'\subset\CM$ as a minimal
$\CU$-submodule containing the lowest weight vector
$[0,\ldots,0]$. The relations ~(\ref{evident'}) show that $\CM'$
is generated from $[0,\ldots,0]$ by the action of operators $e_i,\
i\in I$. Clearly, $\CM'$ is isomorphic to a universal Verma module
over $\CU$.

\begin{conj}
\label{nakos} The class of the structure sheaf $[\CO_{\ul{d}}]$
lies in $\CM'_{\ul{d}}$.
\end{conj}

In what follows we shall assume the validity of Conjecture
~\ref{main'} (as was explained above this is actually not an
assumption for $n>2$).

\begin{prop}
\label{nakosi} The class of $[\D^{-1}_{\ul{d}}]$ lies in
$\CM'_{\ul{d}}$.
\end{prop}

\begin{proof} We have $\CM=\CM'\oplus\CM''$ where $\CM''$ is the
orthogonal complement of $\CM'$ in $\CM$ with respect to the
Shapovalov form. We have to prove that $[\D^{-1}_{\ul{d}}]$ is
orthogonal to $\CM''$. Let $A\in\CM''_{\ul{d}}$. Suppose $A=e_iB$
for some $i\in I$ and $B\in\CM''_{\ul{d}-i}$. Then
$(A,[\D^{-1}_{\ul{d}}])=(B,e_i^*[\D^{-1}_{\ul{d}}])$. Thus up to
(an invertible) monomial in $t,u,v$ we have
$(A,[\D^{-1}_{\ul{d}}])=(B,[D^{-1}_{\ul{d}-i}])$. Hence, arguing
by induction in $\ul{d}$ we may assume that $A\in\CM''_{\ul{d}}$
is orthogonal to the image of any $e_i$. Then $e_i^*A=0$ or,
equivalently, $f_iA=0$ for any $i\in I$. Up to (an invertible)
monomial in $t,u,v$ we have
$(A,[\D^{-1}_{\ul{d}}])=R\Gamma(\CP_{\ul{d}},A)$. Thus we are
reduced to the following claim for $\ul{d}\ne(0,\ldots,0)$:

\begin{equation}
\label{ponjatno} f_iA=0\ \forall\ i\in I\ \Longrightarrow\
R\Gamma(\CP_{\ul{d}},A)=0.
\end{equation}

We will derive ~(\ref{ponjatno}) from the corresponding claim in
the equivariant (complexified) Borel-Moore homology
$H^{\widetilde{T}\times\BC^*\times\BC^*}_{BM}(\CP_{\ul{d}})$. Let
$\on{Td}_{\CP_{\ul{d}}}$ denote the equivariant Todd class in the
completion of the equivariant cohomology. Let also $\on{ch}_*$
denote the homological Chern character map from the equivariant
$K$-theory to the completion of the equivariant Borel-Moore
homology (see e.g. ~\cite{cg}). We define
$$a:=\on{Td}_{\CP_{\ul{d}}}\cup\on{ch}_*A\in
\widehat{H}{}^{\widetilde{T}\times\BC^*\times\BC^*}_{BM}(\CP_{\ul{d}})$$
By the bivariant Riemann-Roch Theorem (see e.g. ~\cite{cg},
~5.11.11) we have
$\on{ch}_*(f_iA)=\on{ch}_*(\bp_*\bq^*A)=\bp_*\bq^*a$ where in the
RHS $\bp_*$ and $\bq^*$ refer to the operations in the (localized
and completed) equivariant Borel-Moore homology. We also have
$R\Gamma(\CP_{\ul{d}},A)=\int_{\CP_{\ul{d}}}a$. Since $\on{ch}_*$
is injective, and the operation
$?\mapsto\on{Td}_{\CP_{\ul{d}}}\cup?$ is invertible, the claim
~(\ref{ponjatno}) follows from the corresponding claim in the
equivariant Borel-Moore homology
$\widehat{H}{}^{\widetilde{T}\times\BC^*\times\BC^*}_{BM}(\CP_{\ul{d}})$:

\begin{equation}
\label{clear} {\mathfrak f}_ia=0\ \forall\ i\in I\
\Longrightarrow\ \int_{\CP_{\ul{d}}}a=0.
\end{equation}

Here ${\mathfrak f}_i=\bp_*\bq^*$ is a part of the action of the
affine Lie algebra $\widehat{\mathfrak{sl}}_n$ on ${\mathfrak
M}:=\oplus_{\ul{d}}
\widehat{\ul{H}}{}^{\widetilde{T}\times\BC^*\times\BC^*}_{BM}(\CP_{\ul{d}})$
(localized and completed equivariant Borel-Moore homology). The
positive generators act as ${\mathfrak e}_i=-\bq_*\bp^*$. This can
be checked along the lines of ~\ref{long}--\ref{longer} but
simpler.

Reversing the argument in the beginning of the proof, we see that
~(\ref{clear}) is equivalent to the statement that the fundamental
cycle $[\CP_{\ul{d}}]\in
\widehat{\ul{H}}{}^{\widetilde{T}\times\BC^*\times\BC^*}_{BM}(\CP_{\ul{d}})$
is contained in the subspace ${\mathfrak M}'$ of ${\mathfrak M}$
generated by the action of ${\mathfrak e}_i,\ i\in I$, from
$[\CP_{(0,\ldots,0)}]$.

Recall the semismall resolution morphism $\pi_{\ul{d}}:\
\CP_{\ul{d}}\to\frP_{\ul{d}}$ to the Uhlenbeck flag space, see
~\cite{fgk}. By the Decomposition Theorem of
Beilinson-Bernstein-Deligne-Gabber, the direct sum of (localized
and completed) equivariant Intersection Homology $'{\mathfrak
M}:=\oplus_{\ul{d}}
\widehat{\ul{IH}}{}^{\widetilde{T}\times\BC^*\times\BC^*}(\frP_{\ul{d}})$
is a direct summand of $\mathfrak M$.

Now ~\cite{b} defines the action of $\widehat{\mathfrak{sl}}_n$ on
$'{\mathfrak M}$, and one can check that the action of ~\cite{b}
is the restriction of the above $\widehat{\mathfrak{sl}}_n$-action
on $\mathfrak M$. It follows that $'{\mathfrak M}={\mathfrak M}'$.
Finally, it is proved in ~\cite{b} that $[\CP_{\ul{d}}]\in\
'{\mathfrak M}$.

This completes the proof of the Proposition.
\end{proof}

\subsection{}
We conclude that $\fu$ is the unique Whittaker vector in the
completion of the Verma module $\CM'$ with the lowest weight
component $\fu_{(0,\ldots,0)}=[(0,\ldots,0)]$ (the common
eigenvector of $e_i^*,\ i\in I$, with the eigenvalue
$(1-v^2)^{-1}$).

Let $\fn'\in\widehat{\CM}{}'$ be the unique common eigenvector of
$f_i,\ i\in I$, with the eigenvalue $(1-v^2)^{-1}$ and with the
lowest weight component $\fn'_{(0,\ldots,0)}=[(0,\ldots,0)]$. Then
$\fn'$ is the orthogonal projection of $\fn$ onto
$\widehat{\CM}{}'$ along $\widehat{\CM}{}''$. Hence the Corollary
~\ref{vot te} yields the following

\begin{cor}
\label{last} One has
$$(\fn'_{\ul{d}},\fu_{\ul{d}})=(-1)^{\sum_{i=0}^{n-1}d_i}
v^{\sum_{i=0}^{n-1}d_i^2-\sum_{i=0}^{n-1}d_id_{i-1}-\sum_{i=0}^{n-1}d_i}
u^{d_0}\prod_{i=1}^nt_i^{d_{i-1}-d_i}[R\Gamma(\CP_{\ul{d}},\CO_{\ul{d}})].$$
\end{cor}

\subsection{Some further remarks}
 The next natural step would be to study the generating function of all
$[R\Gamma(\CP_{\ul{d}},\CO_{\ul{d}})]$'s in a way similar to
subsection ~\ref{generating}; let us denote this function by
$\fJ_{\aff}$. The cohomology (as opposed to $K$-theory) analogue
of this is performed in \cite{b} and \cite{be}. In particular, in
\cite{b} it is shown that such a function is an eigen-function of
a certain linear differential operator of 2nd order (the
``non-stationary analogue" of the quadratic affine Toda
hamiltonian). This fact is used in \cite{be} in order to show that
certain asymptotic of this function is given by the {\it
Seiberg-Witten prepotential} of the corresponding classical affine
Toda system. This agrees well with the results of \cite{neok}
about a similar asymptotic of the partition function of N=2
supersymmetric gauge theory in 4 dimensions.

Unfortunately, in the present ($K$-theoretic) case we can't derive
any good equation for the the function $\fJ_{\aff}$. Thus we do
not know how to generalize the results of \cite{be} to this case.
One can probably show that the results of \cite{neok} on 5d gauge
theory imply that a similar asymptotic (when the classical affine
Toda lattice is replaced by the classical affine relativistic
Toda) is valid for the function $\fJ_{\aff}$, but we do not know
how to derive it from Corollary ~\ref{last}.

\end{document}